\documentclass[a4paper,10pt]{article}
\usepackage{color}

\usepackage{memhfixc}
\usepackage{setspace}

\usepackage[cmex 10]{amsmath}
\usepackage{amsthm}
\usepackage{graphicx}
\usepackage{amsfonts}
\usepackage{amssymb}
\usepackage{afterpage}
\usepackage{natbib}
\usepackage{verbatim} 
\newtheorem{theorem}{Theorem}[section]

\theoremstyle{definition}
\newtheorem{definition}[theorem]{Definition}

\theoremstyle{remark}

\newcommand{\be}{\begin{equation*}}
\newcommand{\ee}{\end{equation*}}
\usepackage[english]{babel}    
\usepackage[dvips]{epsfig}     
\usepackage{epstopdf}
\begin{document}
\title{Efficient Gibbs Sampling for Markov Switching GARCH Models}
\author{\hspace{10pt} Monica Billio \setcounter{footnote}{1}\footnotemark{}\hspace{20pt}
Roberto Casarin\setcounter{footnote}{1}\footnotemark{}\hspace{20pt}Anthony Osuntuyi\setcounter{footnote}{1}\footnotemark{}
\hspace{4pt}\setcounter{footnote}{6}\footnote{Address: Department of Economics, University Ca' Foscari of Venice, Fondamenta San Giobbe 873, 30121, Venice, Italy. Corresponding author: Anthony Osuntuyi, {\rm aosuntuyi@stud.unive.it}. Other contacts: {\rm billio@unive.it} (Monica Billio); {\rm r.casarin@unive.it} (Roberto Casarin).}
\\{\centering \setcounter{footnote}{1}\footnotemark{} University Ca' Foscari of Venice}}
\date{December 2012}
\maketitle
\begin{abstract}
We develop efficient simulation techniques for Bayesian inference on switching GARCH models. Our contribution to existing literature is manifold. First, we discuss different multi-move sampling techniques for Markov Switching (MS) state space models with particular attention to MS-GARCH models. Our multi-move sampling strategy is based on the Forward Filtering Backward Sampling (FFBS) applied to an approximation of MS-GARCH. Another important contribution is the use of multi-point samplers, such as the Multiple-Try Metropolis (MTM) and the Multiple trial Metropolize Independent Sampler, in combination with FFBS for the MS-GARCH process. In this sense we extend to the MS state space models the work of \cite{So06} on efficient MTM sampler for continuous state space models. Finally, we suggest to further improve the sampler efficiency by introducing the antithetic sampling of \cite{crameng2} and \cite{Craiu:2007rr} within the FFBS. Our simulation experiments on MS-GARCH model show that our multi-point and multi-move strategies allow the sampler to gain efficiency when compared with single-move Gibbs sampling.
\end{abstract}
\clearpage        
\section{Introduction}
\indent The study of financial markets volatility has remained a prominent area of research in finance given the important role it plays in a variety of financial problems (e.g. asset pricing and risk management) challenging both investors and fund managers. A remarkable amount of work, ranging from model specification in discrete and continuous time to estimation techniques and finally to applications, have been proposed in the literature. Among volatility models, \cite{Boll86} Generalized Autoregressive Conditional Heteroskedastic (GARCH) model and its variants ranks as the most popular class of models among practitioners. However, from empirical studies, this class of models have been well documented to exhibit high persistence of conditional variance, i.e. the process is close to being non-stationary (nearly integrated). \cite{LamLas90}, among others, argue that the presence of structural changes in the variance process, for which the standard GARCH process cannot account for, may be responsible for this phenomenon. To buttress this point, \cite{MikSta04} estimate a GARCH model on a sample that exhibits structural changes in its conditional variance and obtained a nearly integrated GARCH effect from the estimate. Based on this observation, \cite{HamSum94} and \cite{Cai94} propose a Markov Switching-Autoregressive Conditional Heteroskedastic (MS-ARCH) model, governed by a state variable that follows a first order Markov chain to capture the high volatility persistence, while \cite{Gray96} considers a Markov Switching GARCH (MS-GARCH) model since it can be written as an infinite order ARCH model and may be more parsimonious than the MS-ARCH model for financial data.

The class of MS-GARCH models is gradually becoming a work house among economics and financial practitioners for analysing financial markets data (e.g., see \cite{mar05}). For practical implementation of this class of theoretical models, it is crucial to have reliable parameter estimators. Maximum Likelihood (ML) approach is a natural route to parameter estimation in Econometrics. However, the ML technique is not computationally feasible for MS-GARCH models because of the path dependence problem (see \cite{Gray96}). To this end, \cite{Hen11} and \cite{Bau10} propose Bayesian approach based on Markov Chain Monte Carlo (MCMC) Gibbs technique for estimating the parameters of Markov Switching-Autoregressive Moving Average-Generalized Autoregressive Conditional Heteroskedastic (MS-ARMA-GARCH) and MS-GARCH models respectively. Their proposed algorithm samples each state variable given others individually (single-move Gibbs sampler). This sampler is slowly converging and computationally demanding. Great attention have been paid in the literature at improving such inefficiencies in the context of continuous and possibly non-Gaussian and nonlinear state space models. See, for example, \cite{Fru94}, \cite{Koop00}, \cite{Dejong95} and \cite{cart94} for multi-move Gibbs sampler and \cite{So06} for multi-points and multi-move Gibbs sampling schemes for continuous and nonlinear state space models. To the best of our knowledge there are few works on efficient multi-move sampling scheme for discrete or mixed state space models. See \cite{kim1999state} for a review on multi-move Gibbs for conditionally linear models, \cite{Billio1999} for global Metropolis Hastings algorithm for sampling the hidden states of MS-ARMA models and \cite{Fio12} for multi-move sampling in dynamic mixture models. As regards MS-GARCH models, \cite{Ardia08} develops a Gibbs sampling scheme for the joint sampling of the state variables for the \cite{HMP04} model, which is a particular approximation of a MS-GARCH model, \cite{HeMah10} propose a Sequential Monte Carlo (SMC) algorithm for GARCH models subject to structural breaks, while \cite{Bau11} propose a Particle MCMC (PMCMC) algorithm for estimating GARCH models subject to either structural breaks and regime switching. \cite{Du12}, on the other hand, propose a Metropolis Hastings algorithm for block sampling of the hidden state of infinite state MS-GARCH models. See also \cite{Elliot12} for an alternative approach, i.e. Viterbi-Based technique, for sampling the state variables of MS-GARCH models.

In this paper, we develop an efficient simulation based estimation approach for MS-GARCH models characterized by a finite number of regimes wherein the conditional mean and conditional variance may change over time from one GARCH process to another. We follow a data augmentation framework by including the state variables into the parameter vector. In particular, we propose a Bayesian approach based on MCMC algorithm which allows to circumvent the problem of path dependence by simultaneously generating the states (multi-move Gibbs sampler) from their joint distribution. Our strategy for sampling the state variables is based on Forward Filtering Backward Sampling (FFBS) techniques. As for mixed hidden state models, FFBS algorithm cannot be applied directly on switching GARCH models, we suggest the use of a Metropolis algorithm with an FFBS proposal generated using an auxiliary model. We propose and discuss different auxiliary models obtained by alternative approximations of the MS-GARCH conditional variance equation.

Another original contribution of the paper relates to the Metropolis step for the hidden states. To efficiently estimate MS-GARCH models we consider the class of generalized (multipoint) Metropolis algorithms (see \cite{Liu:2001}, Chapter 5) which extends the standard Metropolis-Hastings (MH) approach (\cite{Hast:70:MCS} and \cite{metrop}). See \cite{Liu:2001} and \cite{RobCas:2007} for an introduction to MH algorithms and a review of various extensions. Multipoint samplers have been proved, both theoretically and computationally, to be effective in improving the mixing rate of the MH chain and the efficiency of the Monte Carlo estimates based on the output of the chain. The main feature of the multipoint samplers is that at each iteration of the MCMC chain the new value of the chain is selected among  multiple proposals, while in the MH algorithm one accepts or rejects a single proposal. In this paper we apply the Multiple-Try Metropolis (MTM) (see \citep{Liu:2000cr}) and some modified MTM algorithms. The superiority of the MTM over standard MH algorithm has been proved in \cite{Craiu:2007rr}, which also propose to apply antithetic and quasi-Monte Carlo techniques to obtain good proposal distributions in the MTM. \cite{So06} applies MTM to the estimation of latent-variable models and finds evidence of superiority of the MTM over standard MH samplers for the latent variable estimation. The author also finds that the efficiency of MTM can further be increased by the use of multi-move sampling. \cite{cascralei:2012} apply the MTM transition to the context of interacting chains. They provide a comparison with standard interacting MH and also estimate the gain of efficiency when using interacting MTM combined with block-sampling for the estimation of stochastic volatility models. We thus combine the MTM sampling strategies with the approximated FFBS techniques for the Markov switching process. In this sense, we extend the work of \cite{So06} to the more complex case of Markov-switching nonlinear state space models. In fact, the use of multiple proposals is particularly suited in this context where the forward filter is used at each iteration to generate only one proposal with a large computational cost. The use of multiple proposals based on the same run of the forward filter is thus discussed. We also apply to this context the antithetic sampling technique proposed by \cite{Craiu:2007rr} to generate correlated proposal within the Multiple-try algorithm, and suggest a Forward Filtering Backward Antithetic Sampling (FFBAS) algorithm which combines the permuted displacement algorithm of \cite{crameng2} with FFBS and possibly produces pairwise negative association among the trajectories of the hidden states. Note that our approach could easily extended to other discrete or mixed state space models.

The paper is organized as follows. Section 2 introduces the MS-GARCH model and discuss inference issues related to existing methods in the literature. In Section 3, we present the Bayesian inference approach and explain the multi-move multipoint sampling strategies. In Section 4, we study the efficiency of our estimation procedure through some simulation experiments. In Section 5, we conclude and discuss possible extensions.  
\section{Markov Switching GARCH models}
\subsection{The model}
A Markov Switching GARCH model is a nonlinear specification of the evolution of a time series assessed to be affected by different states of the world and for which the conditional variance in each state follows a GARCH process. More specifically, let $y_{t}$ be the observed variable (e.g. the return on some financial asset) and $s_{t}$ a discrete, unobserved, state variable which could be interpreted as the state of the world at time $t$. Define $(y_{s},\dots,y_{t})$ and $(s_{s},\dots,s_{t})$ as $y_{s:t}$ and $s_{s:t}$ respectively whenever $s\leq t$ and 0 otherwise. Then
\begin{equation}
y_{t}=\mu_{t}(y_{1:t-1},\theta_{\mu}(s_{t}))+\sigma_{t}(y_{1:t-1},\theta_{\sigma}(s_{t}))\eta_{t},~~~~~~~~~\eta_{t}\overset{iid}{\sim}N(0,1), \label{eq1}
\end{equation}
\begin{equation}
\sigma_{t}^{2}(y_{1:t-1},\theta_{\sigma}(s_{t}))=\gamma(s_{t})+\alpha(s_{t})\epsilon_{t-1}^{2}+\beta(s_{t})\sigma_{t-1}^{2}(y_{1:t-2},\theta_{\sigma}(s_{t-1})),\label{eq2}
\end{equation}
where, $\epsilon_{t}=\sigma_{t}(y_{1:t-1},\theta_{\sigma}(s_{t}))\eta_{t}$, $\theta_{\sigma}(s_{t})=(\gamma(s_{t}),\alpha(s_{t}),\beta(s_{t}))$, $\gamma(s_{t})> 
0$, $\alpha(s_{t})\geq 0$, $\beta(s_{t})\geq 0$, and $s_{t}\in\{1,\ldots,M\}$, $t=1,\dots,T$, is assumed to follow a $M$-state first order Markov chain with transition probabilities $\{\pi_{ij,t}\}_{i,j=1,2,\dots,M}$:
\be
\pi_{ij,t}=p(s_{t}=i|s_{t-1}=j,y_{1:t-1},\theta_{\pi}),~~~\sum_{i=1}^{M}\pi_{ij,t}=1~~~\forall~j=1,2,\dots,M.
\ee
The parameter shift functions $\gamma(s_{t})$, $\alpha(s_{t})$ and $\beta(s_t)$, describe the dependence of parameters 
on the realized regime $s_{t}$ i.e.
\be
\gamma(s_{t})= \sum_{m=1}^{M} \gamma_{m}\mathbb{I}_{s_{t}=m},~~\alpha(s_{t})= \sum_{m=1}^{M} \alpha_{m}\mathbb{I}_{s_{t}=m},~~{\text{and}}~~ \beta(s_{t})= \sum_{m=1}^{M} \beta_{m}\mathbb{I}_{s_{t}=m},
\ee
where,
\be
\mathbb{I}_{s_{t}=m}=
\begin{cases}
1,  &  \text{if}~s_{t}=m\\
0,  &\text{otherwise}\\
\end{cases},
\ee
By defining the allocation variable, $s_{t}$, as a $M$-dimensional discrete vector, $\xi_{t}= (\xi_{1t},\dots,\xi_{Mt})'$, where $\xi_{mt} = \mathbb{I}_{s_{t}=m}$, $m=1,\dots,M,$ the system of equations in (\ref{eq1})-(\ref{eq2}) can be written compactly as
\begin{equation}
y_{t}=\mu_{t}(y_{1:t-1},\xi_{t}'\theta_{\mu})+\sigma_{t}(y_{1:t-1},\xi_{t}'\theta_{\sigma})\eta_{t},~~~~~~~~~\eta_{t}\sim^{iid}N(0,1), \label{eq3}
\end{equation}
\begin{equation}
\sigma_{t}^{2}(y_{1:t-1},\xi_{t}'\theta_{\sigma})=(\xi_{t}'\gamma)+(\xi_{t}'\alpha)\epsilon_{t-1}^{2}+(\xi_{t}'\beta)\sigma_{t-1}^{2}(y_{1:t-2},\xi_{t-1}'\theta_{\sigma}), \label{eq4}
\end{equation}
where $\epsilon_{t}=\sigma_{t}(y_{1:t-1},\xi_{t}'\theta_{\sigma})\eta_{t}$, $\gamma=(\gamma_{1},\dots,\gamma_{M})'$,  $\alpha=(\alpha_{1},\dots,\alpha_{M})'$, $\beta=(\beta_{1},\dots,\beta_{M})'$, $\theta_{\mu}=(\theta_{1\mu},\dots,\theta_{M\mu})'$ and  $\theta_{\sigma}=(\theta_{1\sigma},\dots,\theta_{M\sigma})'$ with $\theta_{m\sigma}=(\gamma_{m},\alpha_{m},\beta_{m})'$ for $m=1,\dots,M$.
for $t=1,\dots,T$. 
Let $\pi_{t}=(\pi_{1t},\dots,\pi_{Mt})$, with $\pi_{it}=(\pi_{i1,t},\dots,\pi_{iM,t})$ for $i=1,2,\dots, M$ and $\sum_{i=1}^{M}\pi_{ij,t}=1$  for all $j=1,2,\dots,M$. Since $\xi_{t}$ follows a $M-$state first order Markov chain, we define the transition probabilities $\{\pi_{ij,t}\}_{i,j=1,2,\dots,M}$ by
\be
\pi_{ij,t}=p(\xi_{t}'=e_{i}'|\xi_{t-1}'=e_{j}',y_{1:t-1},\theta_{\pi}),
\ee
where $e_{i}$ is the $i-$th column of a M-by-M identity matrix. The conditional probability of $\xi_{t}$ given $\xi_{t-1}$, $\theta_{\pi}$ and $y_{1:t-1}$ is given by
\begin{equation}
p(\xi_{t}'|\xi_{t-1}',y_{1:t-1},\theta_{\pi})=\prod_{m=1}^{M}(\pi_{mt}\xi_{t-1})^{\xi_{mt}}, \label{eq5}
\end{equation}
which implies that the probability with which event $m$ occurs at time $t$ is $\pi_{mt}\xi_{t-1}$. 
\subsection{Inference Issues}
Estimating Markov switching GARCH models is a challenging problem since the likelihood of $y_{t}$ depends on the entire sequence of past states up to time $t$ due to the recursive structure of its volatility. To elaborate on this, the likelihood function of the switching GARCH model is given by
\begin{equation}
{\mathcal{L}}(\theta|y_{1:T})\equiv f(y_{1:T}|\theta)=\sum_{i=1}^{M}\dots\sum_{j=1}^{M}f(y_{1:T},\xi_{1}'=e_{i}',\dots,\xi_{T}'=e_{j}'|\theta) \label{eq6}
\end{equation}
where $\theta=(\{\theta_{m\mu},\theta_{m\sigma}\}_{m=1,\dots,M},\theta_{\pi})$. Setting $\xi_{s:t}=(\xi_{s}',\dots,\xi_{t}')$ whenever $s\leq t$, the joint density function of $y_{1:t}$ and $\xi_{1:t}$ on the right hand side of equation (\ref{eq6}) is 
\begin{equation}
\begin{aligned}
f(y_{1:T},\xi_{1:T}|\theta)&=f(y_{1}|\xi_{1:1},\theta_{\mu},\theta_{\sigma})\prod_{t=2}^{T}f(y_{t}|y_{1:t-1},\xi_{1:t},\theta_{\mu},\theta_{\sigma})p(\xi_{t}|y_{1:t-1},\xi_{1:t-1},\theta_{\pi})\\
									&=f(y_{1}|\xi_{1:1},\theta_{\mu},\theta_{\sigma})\prod_{t=2}^{T}f(y_{t}|y_{1:t-1},\xi_{1:t},\theta_{\mu},\theta_{\sigma})\left(\prod_{i=1}^{M}(\pi_{it}\xi_{t-1})^{\xi_{it}}\right),
\end{aligned}\label{eq7}
\end{equation}
with,
\be
\begin{aligned}
f(y_{t}|y_{1:t-1},\xi_{1:t},\theta_{\mu},\theta_{\sigma})&\propto \dfrac{1}{\sigma_{t}(y_{1:t-1},\xi_{t}'\theta_{\sigma})}\exp{\left(-\dfrac{1}{2}\left(\dfrac{y_{t}-\mu_{t}(y_{1:t-1},\xi_{t}'\theta_{\mu})}{\sigma_{t}(y_{1:t-1},\xi_{t}'\theta_{\sigma})}\right)^{2}\right)}.
\end{aligned}
\ee
Given $\sigma_{1}$, recursive substitution in equation (\ref{eq4}) yields
\begin{equation}
\sigma_{t}^{2}=\sum_{i=0}^{t-2}\left[\xi_{t-i}'\gamma+(\xi_{t-i}'\alpha)\epsilon_{t-1-i}^{2}\right]\prod_{j=0}^{i-1}\xi_{t-j}'\beta+\sigma_{1}^{2}\prod_{i=0}^{t-2}\xi_{t-i}'\beta. \label{eq8}
\end{equation}
Equation (\ref{eq8}) clearly shows the dependence of conditional variance at time $t$ on the entire history of the regimes and by inference the dependence of the likelihood function on the entire history of the regimes. The evaluation of the likelihood function over a sample of length $T$, as can be seen in equation (\ref{eq6}), involves integration (summation) over all $M^{T}$ unobserved states i.e. integration over all $M^{T}$ possible (unobserved) regime paths. This requirement makes the maximum likelihood estimation of equation (\ref{eq6}) infeasible in practice. 

Two major approaches have been developed in the literature in order to circumvent this path dependence problem. One approach involves the use of model approximation while the other is simulation based. 

As regards to the model approximation approach, \cite{Cai94} and \cite{HamSum94} approximated the MS-GARCH model by an MS-ARCH model. This approach effectively makes the model tractable because the lagged conditional variance that makes the conditional variance dependent on the history of regime has been dropped. \cite{KauFrau02} employed the algorithm developed by \cite{Chib96} for a Markov mixture models to compute the marginal likelihood of the MS-ARCH model but noted that this methodology cannot be carried over to the MS-GARCH model because of the path dependence problem. Another approximation approach can be credited to \cite{Gray96} who noted that the conditional density of the return is essentially a mixture of distributions with time-varying mixing parameter and in particular under normality assumption he suggested the use of aggregate conditional variances over all regimes as the lagged conditional variance when constructing the conditional variance at each time step. Extensions of \cite{Gray96} model can be found in \cite{Due97}, \cite{klas02} and  \cite{HMP04} among others. \cite{Abracoh07} provide stationarity conditions for some of these approximations. The problem with this approach is that these approximations cannot be verified.

Among the simulation based approaches proposed in the literature there is the Bayesian estimation technique by \cite{Bau10}. In particular, they develop a single-move MCMC Gibbs sampler for a Markov switching GARCH model with a fixed number of regimes. The authors also provide sufficient conditions for geometric ergodicity and existence of moments of the process. Their estimation approach, though quite promising, has one main limitation that has rendered it unattractive. The single-move Gibbs sampler is inefficient i.e. draws from the single-move scheme are noted to be highly correlated and thus slow down the convergence of the Markov chain. An alternative simulation based approach is the particle filter approach proposed by \cite{HeMah10}. They develop a sequential Monte Carlo method for estimating GARCH models subject to an unknown number of structural breaks. 

In the next section, we propose an efficient Bayesian estimation procedure for estimating the parameters of MS-GARCH models by simultaneously generating the whole state vector.
\section{Bayesian Inference}
Based on the aforementioned inference issues associated with MS-GARCH models, we present a Bayesian approach based on MCMC Gibbs algorithm which allows us to circumvent the path dependence problem and efficiently sample the state trajectory. The purpose of this algorithm is to generate samples from the posterior distribution which are then used for its characterization. We follow a data augmentation framework by treating the state variables as parameters of the model and construct the likelihood function assuming the states known. 

Before proceeding with the elicitation of our proposed Bayesian technique, it is important that we make explicit the parametric specification of the conditional mean, $\mu_{t}(y_{1:t-1},\xi_{t}'\theta_{\mu})$, of the return process $y_{t}$ in equation (\ref{eq3}) and the transition probabilities $p(\xi_{t}'|\xi_{t-1}',y_{1:t-1},\theta_{\pi})$. Since our major aim is to define a technique for sampling the state variables efficiently, which in turn will affect other parameter estimates, we assume for expository purposes a conditional mean defined by a constant switching parameter given by $\xi_{t}'{\mu}$  where $\mu=(\mu_{1},\dots,\mu_{M})'$ and constant transition probabilities. Alternative specification such as switching ARMA process could be thought of for the conditional mean and time varying transition probabilities may be defined by following \cite{Gray96} approach, i.e. specifying transition probabilities as a function of past observables. Under this specification, the augmented parameter set of our model consists of $\xi_{1:T}$, $\theta = (\theta_{\mu},\theta_{\sigma},\theta_{\pi})$ where $\theta_{\mu}=\mu$, $\theta_{\pi}=(\{\pi_{m}\}_{m=1,\dots,M})$ and $\theta_{\sigma}=(\{\theta_{m\sigma}\}_{m=1,\dots,M})$ with $\theta_{m\sigma}=(\gamma_{m},\alpha_{m},\beta_{m})$, $\pi_{m}=(\pi_{m1},\dots,\pi_{mM})$ and $\sum_{m=1}^{M}\pi_{mm^{*}}=1$ $\forall$ $m^{*}=1,\dots,M$. The prior distributions of the parameter vector are assumed to be independent and chosen as follows
\be
\begin{aligned}
\theta_{\pi}&\sim  \prod_{m=1}^{M}{\text{Dirichlet}}(\nu_{1m},\dots,\nu_{Mm})\\
\theta_{\mu}&\sim \prod_{m=1}^{M} \mathcal{U}_{[a_{m\mu},b_{m\mu}]}\\
\theta_{\sigma}&\sim \prod_{m=1}^{M} \mathcal{U}_{[a_{m\gamma},b_{m\gamma}]}\mathcal{U}_{[a_{m\alpha},b_{m\alpha}]}\mathcal{U}_{[a_{m\beta},b_{m\beta}]}.
\end{aligned}
\ee
where $\nu_{1m},\dots,\nu_{Mm},a_{m\mu},b_{m\mu},a_{m\gamma},b_{m\gamma},a_{m\alpha},b_{m\alpha},a_{m\beta},b_{m\beta}$ $\forall$ $m=1,\dots,M$ are hyperparameters to be defined. The supports of the prior distribution of $\theta_{\mu}$ and $\theta_{\sigma}$ will be chosen to avoid label switching (identifiability restriction). See \cite{Fru06} for an introduction to label switching problem for dynamic mixtures and MS models and\cite{Bau10} for illustration of the identification constraint for MS-GARCH models. The choice of the prior supports also helps in preventing regime degeneration. The joint prior distribution is thus proportional to 
\begin{equation}
f(\theta) \propto \prod_{m=1}^{M}{\text{Dirichlet}}(\nu_{1m},\dots,\nu_{Mm}).
\end{equation}
The posterior density of the augmented parameter vector given by 
\begin{equation}
\begin{aligned}
f(\theta,\xi_{1:T}|y_{1:T})&\propto f(y_{1:T},\xi_{1:T},\theta)\\
												 &= f(y_{1:T}|\xi_{1:T},\theta)f(\xi_{1:T}|\theta)f(\theta).\label{eq9a}
\end{aligned}
\end{equation}
cannot be identified with any standard distribution, hence we cannot sample directly from it. Using Gibbs sampler, we can generate samples from this high-dimensional posterior density. This will be done by iteratively sampling from the following three full conditional distributions
\begin{itemize}
\item $p(\xi_{1:T}|\theta,y_{1:T})$,
\item $f(\theta_{\pi}|\theta_{\mu},\theta_{\sigma},\xi_{1:T},y_{1:T})= f(\theta_{\pi}|\xi_{1:T})$, and
\item $f(\theta_{\sigma},\theta_{\mu}|\theta_{\pi},\xi_{1:T},y_{1:T})= f(\theta_{\sigma},\theta_{\mu}|\xi_{1:T},y_{1:T})$.
\end{itemize}
These full conditional distributions are easier to manage and sample from because they can either be associated with a known distribution or simulated by a lower dimensional auxiliary sampler. In the following subsections we present in details our sampling procedure.
\subsection{Sampling the state variables $\xi_{1:T}$.}
To sample $\xi_{1:T}$ using the single move algorithm, one relies on computing
\begin{equation}
p(\xi_{t}|\xi_{1:t-1},\xi_{t+1:T},\theta,y_{1:T})\propto \prod_{m=1}^{M} \left(\pi_{m}\xi_{t-1}\right)^{\xi_{mt}}\left(\pi_{m}\xi_{t}\right)^{\xi_{m,t+1}}\prod_{j=t}^{T}f(y_{j}|\xi_{j},\theta,y_{1:j-1})\label{eqsg}
\end{equation}
for each value $\xi_{t}$ in $\{e_{m}:m=1,\dots,M\}$ and dividing each evaluation by the sum of the $M$ points to get the normalized discrete distribution of $\xi_{t}$ from which to sample. Sampling from such a distribution once the probabilities are known is similar to sampling from a Multinomial distribution.
On the other hand, the full joint conditional distribution of the state variables, $\xi_{1:T}$, given the parameter values and return series
\begin{equation}
p(\xi_{1:T}|\theta,y_{1:T}) \propto f(y_{1:T}|\xi_{1:T},\theta)p(\xi_{1:T}|\theta)
\end{equation}
is a non-standard distribution. Therefore multi-move sampling is not feasible.  For this reason, we consider a generalization of MH (i.e. multipoint Metropolis-Hastings) strategy for generating proposals for the state variables. Multipoint samplers are designed to consider multiple proposals at each iteration of an MH and to choose the new value of the chain from this trial set. The multi-move and multipoint sampling procedures are of interest because of their potentials at addressing issues associated with multi-modality of the target function (i.e. in the event that the target distribution is multi-modal in nature the MCMC chain runs the risk of getting trapped in local modes) and autocorrelation of samples from the Metropolis-Hasting's chain. Our scheme generally involves running a FFBS on the auxiliary sampler to generate several proposals at each iteration step. Let the proposal distribution be denoted by 
\begin{equation}
q(\xi_{1:T}|\theta,y_{1:T})=q(\xi_{T}'|\theta,y_{1:T})\prod_{t=1}^{T-1}q(\xi_{t}'|\xi_{t+1}',\theta,y_{1:t}),\label{eq9}
\end{equation}
where $q(\xi_{t}'|\xi_{t+1}',\theta,y_{1:t})\propto q(\xi_{t}'|y_{1:t},\theta)q(x_{t+1}'|x_{t}',\theta)$ with $q(\xi_{t}'|y_{1:t},\theta)$ representing filtered probability. A discussion on the proposal distribution is presented in section 3.2. In the following, we discuss the three multipoint algorithms considered in this paper. 
\subsubsection{Multiple-Try Metropolis Sampler}
\cite{Liu:2000cr} suggest the Multiple-Try Metropolis (MTM) sampler scheme. As in the general case of multipoint samplers, their idea is to consider several points generated by a proposal distribution so that possibly a larger region from which the new value for the chain is chosen can be investigated. By using the multiple-try strategy, it is easier for the iterates to jump from one local maximum to another and thus speed up the convergence to the desired target distribution. Samples from the proposal distribution will be generated by FFBS algorithm. We present below a sketch of the main ingredients needed in Forward Filter (FF) and Backward Sampling (BS) algorithm and refer the reader to \cite{Fru06} for detailed presentation of this procedure. At time $t$, given $\theta$ and $y_{1:t}$ the FF probabilities are obtained by first computing the one-step ahead prediction 
\be
q(\xi_{t}'|,\theta,y_{1:t-1})=\sum_{i=1}^{M}\left(\prod_{j=1}^{M}(\pi_{j}e_{i})^{\xi_{j,t}}\right)q(\xi_{t-1}'=e_{i}'|\theta,y_{1:t-1}),
\ee
then, the FF is
\begin{equation}
q(\xi_{t}'|\theta,y_{1:t})=\dfrac{g(y_{t}|\xi_{t}',\theta,y_{1:t-1})q(\xi_{t}'|\theta,y_{1:t-1})}{\sum_{i=1}^{M}g(y_{t}|\xi_{t}'=e_{i}',\theta,y_{1:t-1})q(\xi_{t}'=e_{i}'|\theta,y_{1:t-1})},
\end{equation}
where $g(y_{t}|\xi_{t}',\theta,y_{1:t-1})$ is the conditional density of the return process under the auxiliary model. Using the output of the FF, we compute  $q(\xi_{T}'|\theta,y_{1:T})$ and
\begin{equation}
q(\xi_{t}'|\xi_{t+1}',\theta,y_{1:t})=\dfrac{\prod_{j=1}^{M}\left(\pi_{j}\xi_{t}\right)^{\xi_{j,t+1}} q(\xi_{t}'|\theta,y_{1:t})}{q(\xi_{t+1}'|\theta,y_{1:t})},\label{eqBS}
\end{equation} 
for $t=T-1,T-2,\dots,2,1$. Then at each time step we sample $\xi_{T}'$ from $q(\xi_{T}'|\theta,y_{1:T})$ and $\xi_{t}'$ from $q(\xi_{t}'|\xi_{t+1}',\theta,y_{1:t})$ iteratively for $t=T-1,T-2,\dots,2,1$. This is the BS step. The BS procedure is implemented by first noting that $\xi_{t+1}$ is the most recent value sampled for the hidden Markov chain at $t+1$ and since $\xi_{t}$ can take one of $e_{1},\dots,e_{M}$, we compute the expression in equation (\ref{eqBS}) for each of these values. Then sampling $\xi_{t}'$ from $q(\xi_{t}'|\xi_{t+1}',\theta,y_{1:t})$ once the corresponding probabilities for $\xi_{i}'=e_{i}'$ for $i=1,\dots,M$ are known may be compared to sampling from a multinomial distribution. Note that at each iteration step of the MCMC procedure we only need a single run of the Forward Filter (FF) for generating generating multiple proposals using Backward Sampling (BS).  

A summary of our MTM algorithm is given in algorithm 1.

		\noindent\hrulefill\\
		{\bf Algorithm 1} MTM Sampler
   	\begin{itemize}
			\item[i.] Choose a starting value $\xi_{1:T}^{0}$.
			\item[ii.] Let $\xi_{1:T}^{(r-1)}$ be the value of the MTM at the $(r-1)$-th iteration.
			\item[iii.] Construct a trial set $\{\xi_{1:T,1},\xi_{1:T,2},\dots,\xi_{1:T,K}\}$ containing $K$ state variable paths drawn from the proposal distribution $q(\xi_{1:T}|\theta^{(r-1)},y_{1:T})$. 
			\item[iv.] Evaluate 
					\be
					W_{k}(\xi_{1:T,k},\xi_{1:T}^{(r-1)})=\dfrac{p(\xi_{1:T,k}|\theta^{(r-1)},y_{1:T})}{q(\xi_{1:T,k}|\theta^{(r-1)},y_{1:T})},~~\forall 						k=1,\dots, K.
					\ee
			\item[v.] Select ${\tilde{\xi}}_{1:T}$ from $\{\xi_{1:T,1},\xi_{1:T,2}\dots,\xi_{1:T,K}\}$ according to the probability
					\be
					p_{k}=\dfrac{W_{k}(\xi_{1:T,k},\xi_{1:T}^{(r-1)})}{\sum_{k=1}^{K}W_{k}(\xi_{1:T,k},\xi_{1:T}^{(r-1)})},~~\forall k=1,\dots, K.
					\ee
			\item[vi.] Construct a reference set $\{\xi_{1:T,1}^{*},\xi_{1:T,2}^{*},\dots,\xi_{1:T,K}^{*}\}$ by setting the first $K-1$ elements to a new set of samples drawn from the proposal distribution $q(\xi_{1:T}|\theta^{(r-1)},y_{1:T})$ and the $K-$th element $\xi_{1:T,K}^{*}$ to $\xi_{1:T}^{(r-1)}$.
			\item[vii.] Draw $u\sim\mathcal{U}_{[0,1]}$.
			\item[viii.] Set 
					\be
					\xi_{1:T}^{(r)} =
					\begin{cases}
 	 				{\tilde{\xi}}_{1:T} & \text{if}~~u\le \alpha({\tilde{\xi}}_{1:T},\xi_{1:T}^{(r-1)})\\
 					\xi_{1:T}^{(r-1)} & \text{otherwise}
					\end{cases}
					\ee
					where,
					\be
		\alpha({\tilde{\xi}}_{1:T},\xi_{1:T}^{(r-1)})=\min\left(1,\dfrac{\sum_{k=1}^{K}W_{k}(\xi_{1:T,k},\xi_{1:T}^{(r-1)})}{\sum_{k=1}^{K}W_{k}(\xi_{1:T,k}^{*},{\tilde{\xi}}_{1:T})}\right).
					\ee				
		\end{itemize}
   	\noindent\hrulefill\\

Observe that the MTM algorithm reduces to standard Metropolis-Hasting algorithm when $K=1$. We also note that alternative weight function other than the importance weight function assumed in the MTM algorithm presented above could be defined.
\subsubsection{Multiple-trial Metropolized Independent Sampler (MTMIS)}
As we are using independent proposal distributions in the MTM algorithm, the generation of the set of reference points is not needed to have a possibly more efficient generalized MH algorithm. Thus, following the suggestion of \cite{Liu:2001} we combine the MTM with the metropolized indpendent sampler and obtain Algorithm 2. The main advantage is that one can use multiple proposals without generating the reference points, obtaining thus a decrease of the computational complexity of the algorithm.

    \noindent\hrulefill\\
		{\bf Algorithm 2} Multiple-trial Metropolized independent Sampler (MTMIS)
   	\begin{itemize}
			\item[i.] Choose a starting value $\xi_{1:T}^{0}$.
			\item[ii.] Let $\xi_{1:T}^{(r-1)}$ be the value of the MTM at the $(r-1)$-th iteration.
			\item[iii.] Construct a trial set $\{\xi_{1:T,1},\xi_{1:T,2},\dots,\xi_{1:T,K}\}$ containing $K$ state variable paths drawn from the proposal distribution. 			\item[iv.] Evaluate 
					\be
					W_{k}(\xi_{1:T,k})=\dfrac{p(\xi_{1:T,k}|,\theta^{(r-1)},y_{1:T})}{q(\xi_{1:T,k}|\theta^{(r-1)},y_{1:T})},~~\forall~~k=1,\dots, K,~~~\text{and define}~~~W = \sum_{k=1}^{K}W_{k}(\xi_{1:T,k})
					\ee					
			\item[v.] Select ${\tilde{\xi}}_{1:T}$ from $\{\xi_{1:T,1},\xi_{1:T,2}\dots,\xi_{1:T,K}\}$ according to the probability
					\be
					p_{k}=\dfrac{W_{k}(\xi_{1:T,k})}{\sum_{k=1}^{K}W_{k}(\xi_{1:T,k})},~~\forall k=1,\dots, K.
					\ee
			\item[vi.] Draw $u\sim\mathcal{U}_{[0,1]}$.
			\item[vii.] Set 
					\be
					\xi_{1:T}^{(r)} =
					\begin{cases}
 	 				{\tilde{\xi}}_{1:T} & \text{if}~~u\le \alpha({\tilde{\xi}}_{1:T},\xi_{1:T}^{(r-1)})\\
 					\xi_{1:T}^{(r-1)} & \text{otherwise}
					\end{cases}
					\ee
					where,
					\be
		\alpha({\tilde{\xi}}_{1:T},\xi_{1:T}^{(r-1)})=\min\left(1,\dfrac{W}{W-W({\tilde{\xi}}_{1:T})+W({\xi}_{1:T}^{(r-1)})}\right).
					\ee				
		\end{itemize}
   	\noindent\hrulefill\\
\subsubsection{Multiple Correlated-Try Metropolis (MCTM) Sampler}
To further improve the efficiency the MTM algorithm and to ensure that a larger portion of the sample space is explored for better mixing and shorter running time, we propose the use of correlated proposals. There are various ways of introducing correlation among proposals e.g. antithetic and stratified approaches. In this paper, we study the antithetic approach. The use of antithetic sampling in a Gibbs sampling context allows for a gain of efficiency. \cite{pitt96} propose a blocking method with antithetic approach for non-Gaussian state space models, \cite{Jasra09} propose a scheme for reducing the variance of estimates from the standard Metropolis-within-Gibbs sampler by introducing antithetic samples while \cite{Olsson08} propose a forward filtering backward smoothing particle filter algorithm with antithetic proposal. Here we follow \cite{Craiu:2007rr} which use antithetic proposals within a multi-point sampler and apply their idea to the context of discrete state space models. We propose a correlated proposal MTM sampler based on a combination of the FFBS and antithetic sampling techniques. To the best of our knowledge, antithetic proposals of this kind have not been used in the context of Markov switching nonlinear state space models. The idea is to choose, at each step of the MCMC algorithm, a new hidden state trajectory from negatively correlated proposals instead of independent proposals. Following the suggestion of \cite{Liu:2001}, we obtain Algorithm 3.

 \noindent\hrulefill\\
		{\bf Algorithm 3} Multiple Correlated-Try Metropolis (MCTM) Sampler 
   	\begin{itemize}
			\item[i.] Choose a starting value $\xi_{1:T}^{0}$.
			\item[ii.] Let $\xi_{1:T}^{(r-1)}$ be the value of the MTM at the $(r-1)$-th iteration.
			\item[iii.] Construct a trial set $\{\xi_{1:T,1},\xi_{1:T,2},\dots,\xi_{1:T,K}\}$ containing $K$ correlated state variable paths drawn from the proposal distribution. 
			\item[iv.] Evaluate 
					\be
					\begin{aligned}
					W_{1}(\xi_{1:T,1})&=\dfrac{p(\xi_{1:T,1}|\theta^{(r-1)},y_{1:T})}{q(\xi_{1:T,1}|\theta^{(r-1)},y_{1:T})},\\
					W_{k}(\xi_{1:T,1:k})&= W_{k-1}(\xi_{1:T,1:k-1})\dfrac{p(\xi_{1:T,k-1}|\theta^{(r-1)},y_{1:T})}{q(\xi_{1:T,k-1}|\theta^{(r-1)},y_{1:T})}~~\forall~~k=2,\dots, K,
					\end{aligned}
					\ee					
			\item[v.] Select ${\tilde{\xi}}_{1:T}$ from $\{\xi_{1:T,1},\xi_{1:T,2}\dots,\xi_{1:T,K}\}$ according to the probability
					\be
					p_{k}=\dfrac{W_{k}(\xi_{1:T,1:k},\xi_{1:T}^{(r-1)})}{\sum_{k=1}^{K}W_{k}(\xi_{1:T,1:k},\xi_{1:T}^{(r-1)})},~~\forall k=1,\dots, K.
					\ee
			\item[vi.] Supposing ${\tilde{\xi}}_{1:T}=\xi_{1:T,l}$ is chosen in item (v) above, create a reference set\\ $\{\xi_{1:T,1}^{*},\xi_{1:T,2}^{*},\dots,\xi_{1:T,K}^{*}\}$ by letting 
			\be
			\begin{aligned}
			\xi_{1:T,j}^{*} &=\xi_{1:T,l-1} ~~~\forall~~ j=1,\dots,l-1\\
			\xi_{1:T,l}^{*} &=\xi_{1:T}^{(r-1)}			
			\end{aligned}			
			\ee
			and drawing $\xi_{1:T,j}^{*}$ for $j=l+1,\dots,K$ from the proposal distribution.			
			\item[vii.] Draw $u\sim\mathcal{U}_{[0,1]}$.
			\item[viii.] Set 
					\be
					\xi_{1:T}^{(r)} =
					\begin{cases}
 	 				{\tilde{\xi}}_{1:T} & \text{if}~~u\le \alpha({\tilde{\xi}}_{1:T},\xi_{1:T}^{(r-1)})\\
 					\xi_{1:T}^{(r-1)} & \text{otherwise}
					\end{cases}
					\ee
					where,
					\be
		\alpha({\tilde{\xi}}_{1:T},\xi_{1:T}^{(r-1)})=\min\left(1,\dfrac{\sum_{k=1}^{K}W_{k}(\xi_{1:T,1:k},\xi_{1:T}^{(r-1)})}{\sum_{k=1}^{K}W_{k}(\xi_{1:T,1:k}^{*},{\tilde{\xi}}_{1:T})}\right).
					\ee					
		\end{itemize}
   	\noindent\hrulefill\\

The simplest way to introduce negative correlation between the trajectories generated with the FFBS algorithm is to use, at a given iteration $r$ of the sampler and for the $t$-th hidden state, a set of $K$ uniform random numbers $U_{t,k}^{(r)}$, $k=1,\ldots,K$ generated following the permuted displacement method (see \cite{ArvJoh82} and \cite{crameng2}) given in Algorithm 4. The uniform random numbers are then use within the BS procedure to generate correlated proposals.

		\noindent\hrulefill\\
		{\bf Algorithm 4} Permuted displacement method
		\begin{itemize}
		\item Draw $r_{1}\sim\mathcal{U}_{[0,1]}$
		\item For $k=2,\ldots,K-1$, set $r_{k}=\lfloor 2^{k-2}r_{1}+1/2\rfloor$ where $\lfloor x \rfloor$ denotes the fractional part of $x$
		\item Set $r_{K}=1-\{2^{K-2}r_{1}\}$
		\item Pick at random $\sigma\in S_{K}$, where $S_{K}$ is the set of all possible permutation of the integers $\{1,\ldots,K\}$
		\item For $k=1,\ldots,K$, set $U_{k}=r_{\sigma(k)}$
		\end{itemize}
		\noindent\hrulefill

For $K=3$, \cite{crameng2} show that the random numbers generated with the permuted displacement method are pairwise negatively associated (PNA). The definition of PNA given in the following is adopted from \cite{crameng2}.
\begin{definition}[pairwise negative association]
The random variables $\xi_{t,1}$,$\xi_{t,2}$.\dots,$\xi_{t,K}$ are said to be pairwise negatively associated (PNA) if, for any nondecreasing functions $f_{1}$, $f_{2}$ and $(i,j)\in\{1,\dots,K\}^{2}$ such that $i\neq j$
\be
Cov(f_{1}(\xi_{t,i}),f_{2}(\xi_{t,j}))\leq 0
\ee
whenever this covariance is well defined.
\end{definition}
The proof for the case $K\geq4$ is still an open issue. For this reason we consider in our algorithm $K\leq3$. The presence of PNA in the case of $K\geq4$ proposals depends on the degrees of uniformity of the filtering probability and the gain of efficiency should be proved computationally in each applications.

We use the permuted sampler to generate $K=2$ multi-move and correlated proposals in the backward sampling step of the FFBS. In order to show how the antithetic sampler works, we consider the case where the hidden Markov switching process has two states, i.e. $\xi_{t}=(\xi_{1t},\xi_{2t})'$ and for notational convenience let $\{q_{t}^{(r)}\}_{t=1:T}$ be the sequence of filtered probabilities of being in state 1 at the $r$-th iteration of the sampler, then we define the backward antithetic samples $\xi_{t,1}$ and $\xi_{t,2}$ as follows
\be
\xi_{t,1}=
\begin{pmatrix}
\mathbb{I}_{U_{t}^{(r)}<q_{t}^{(r)}}\\
\mathbb{I}_{U_{t}^{(r)}\geq q_{t}^{(r)}}\\
\end{pmatrix},
\quad\quad\quad
\xi_{t,2}=
\begin{pmatrix}
\mathbb{I}_{V_{t}^{(r)}<q_{t}^{(r)}}\\
\mathbb{I}_{V_{t}^{(r)}\geq q_{t}^{(r)}}\\
\end{pmatrix}
\ee
where $V_{t}^{(r)}=1-U_{t}^{(r)}$. It is possible to show that
\be
Cov(\xi_{t,1}^{(r)},\xi_{t,2}^{(r)})=
\begin{pmatrix}
(2q_{t}^{(r)}-1)\mathbb{I}_{q_{t}^{(r)}>\frac{1}{2}}-\left(q_{t}^{(r)}\right)^{2} & \left(q_{t}^{(r)}(1-q_{t}^{(r)})\right)^{2}\\
\left(q_{t}^{(r)}(1-q_{t}^{(r)})\right)^{2}  & (1-2q_{t}^{(r)})\mathbb{I}_{q_{t}^{(r)}<\frac{1}{2}}-\left(1-q_{t}^{(r)}\right)^{2}
\end{pmatrix}.
\ee
Using the expected value of the square of the Euclidean distance, $d(\xi_{t,1},\xi_{t,2})$, between this two antithetic samples to investigate the nature of the  antithetic samples, extremely antithetic proposals is obtained when the distance on average is optimal.  
\begin{equation}
E[d^{2}(\xi_{t,1},\xi_{t,2})]= 2 -2\left((2q_{t}^{(r)}-1)\mathbb{I}_{q_{t}^{(r)}>\frac{1}{2}}  + (1-2q_{t}^{(r)})\mathbb{I}_{q_{t}^{(r)}<\frac{1}{2}}\right).\label{eqaa}
\end{equation}
From equation (\ref{eqaa}) extreme antithetic is attained when $q_{t}^{(r)}$ is equal to 0.5, which can be easily found in applications where regimes
exhibit similar persistence level.. 
\subsection{Auxiliary models for defining the proposal distribution}
In order to built proposal distributions for the state variables, we will exploit all the knowledge we have about the full conditional distribution. The first step is to approximate the MS-GARCH model by eliminating the problem of path dependence and then deriving a proposal distribution for state variables from the auxiliary model does obtained. A possible way of circumventing the path dependence problem inherent in the MS-GARCH model is to replace the lagged conditional variance appearing in the definition of the GARCH model with a proxy. A look into the literature shows different auxiliary models which differs only by the content of the information used in defining the proxy used in each case. In general, various MS-GARCH (as available in the literature) can be obtained by approximating the conditional variance
\be
\sigma_{t}^{2}(y_{1:t-1},\theta_{\sigma}(s_{t}))= V[y_{t}|y_{1:t-1},s_{1:t}] = V[\epsilon_{t}|y_{1:t-1},s_{1:t}]
\ee
of the GARCH process as follow
\begin{equation}
{\sigma_{t}}^{2}(y_{1:t-1},\xi_{t}'\theta_{\sigma})\approx\xi_{t}'\gamma + (\xi_{t}'\alpha)\epsilon_{(X)t-1}^{2} + (\xi_{t}'\beta)\sigma_{(X)t-1}^{2}. \label{eq10}
\end{equation}
In the subsection we present alternative specifications of $\epsilon_{(X)t-1}$ and $\sigma_{(X)t-1}^{2}$ that define different approximations of the MS-GARCH model. The variable $X$ can take on any of $B,G,D,SK,K$ with each notation representing, respectively, the Basic approximation, \cite{Gray96} approximation, \cite{Due97} approximation, Simplified version of \cite{klas02} approximation and \cite{klas02} approximation.
\subsubsection{Model 1}
As a first attempt at eliminating the path dependent problem, we note that the conditional density of $\epsilon_{t}$ is  a mixture of normal distribution with zero mean and time varying variance. Hence, we approximate the switching GARCH model by replacing the lagged conditional variance, $\sigma^{2}_{t-1}$, with the variance $\sigma^{2}_{(B)t-1}$ of the conditional density of $\epsilon_{t}$ i.e.
\be
\begin{aligned}
\epsilon_{(B)t-1}&= y_{t-1}-\mu_{(B)t-1}\\
\mu_{(B)t-1}&=E[\mu_{t-1}(y_{1:t-2},\xi_{t-1}'\theta_{\mu})|y_{1:t-2}] = E[y_{t-1}|y_{1:t-2}]\\
            &=\sum_{m=1}^{M}\mu_{t-1}(y_{1:t-2},e_{m}'\theta_{\mu}) q(\xi_{t-1}'=e_{m}'|y_{1:t-2}),\\
\sigma_{(B)t-1}^{2}&=E[\sigma_{t-1}^{2}(y_{1:t-2},\xi_{t-1}'\theta_{\sigma})|y_{1:t-2}]=E[\epsilon_{t-1}^{2}|y_{1:t-2}]=V(\epsilon_{t-1}|y_{1:t-2})\\
							&=\sum_{m=1}^{M}\sigma_{t-1}^{2}(y_{1:t-2},e_{m}'\theta_{\sigma})q(\xi_{t-1}'=e_{m}'|y_{1:t-2}).
\end{aligned}
\ee 
Observe that in this approximation scheme $\mu_{(B)t-1}$ and $\sigma_{(B)t-1}^{2}$ are functions of $y_{1:t-2}$ and the information coming from $y_{t-1}$ is lost. With $q(\xi_{t-1}'=e_{m}'|y_{1:t-2})$ known for $m=1,\dots,M$, $\mu_{(B)t-1}$ can easily be computed while $\sigma_{(B)t-1}^{2}$ can be computed recursively since  ${\sigma}_{t-1}^{2}(y_{1:t-2},e_{m}'\theta_{\sigma})$ depends on $\sigma_{(B)t-2}^{2}$. Note that in this approximation the conditioning is on $y_{1:t-2}$. This approach represents a starting point for other approximations hence we tag it Basic Approximation.
\subsubsection{Model 2}
\cite{Gray96} notes that the conditional density of the return process, $y_{t}$, of the switching GARCH model is a mixture of normal distribution with time-varying parameters. Hence, he suggests the use of the variance of the conditional density $\sigma^{2}_{(G)t-1}$ of $y_{t}$ as a proxy for the lagged of the conditional variance $\sigma^{2}_{t-1}$ switching GARCH process i.e.
\be
\begin{aligned}
\epsilon_{(G)t-1}&=y_{t-1}-\mu_{(G)t-1}\\
\mu_{(G)t-1}&=\mu_{(B)t-1}\\
\sigma_{(G)t-1}^{2}&=V(y_{t-1}|y_{1:t-2})=V\left(E[y_{t-1}|y_{1:t-2},\xi_{t-1}']|y_{1:t-2}\right)+E[V\left(y_{t-1}|y_{1:t-2},\xi_{t-1}'\right)|y_{1:t-2}]\\
									 &= V(\mu_{t-1}(y_{1:t-2},\xi_{t-1}'\theta_{\mu})|y_{1:t-2}) + E[\sigma_{t-1}^{2}(y_{1:t-2},\xi_{t-1}'\theta_{\sigma})|y_{1:t-2}]\\
									 &=E[(\mu_{t-1}(y_{1:t-2},\xi_{t-1}'\theta_{\mu}))^{2}|y_{1:t-2}] - (E[\mu_{t-1}(y_{1:t-2},\xi_{t-1}'\theta_{\mu})|y_{1:t-2}])^{2} + \sigma_{(B)t-1}^{2}\\
									 &=\sum_{m=1}^{M}(\mu_{t-1}(y_{1:t-2},e_{m}'\theta_{\mu}))^{2}q(\xi_{t-1}'=e_{m}'|y_{1:t-2}) - (\mu_{(B)t-1})^{2} + \sigma_{(B)t-1}^{2}.
\end{aligned}
\ee 
Similarly, as in model 1, information on $y_{t-1}$ is lost in this approximation scheme as $\mu_{(G)t-1}$ and $\sigma_{(G)t-1}^{2}$ are functions of $y_{1:t-2}$. By recursion, $\sigma_{(G)t-1}^{2}$ can be computed since $\sigma_{(B)t-1}^{2}$ depends on $\sigma_{(G)t-2}^{2}$ through $\sigma_{t-1}^{2}(y_{1:t-2},e_{m}'\theta_{\sigma})$. Within this framework the conditioning is also on $y_{1:t-2}$. The major difference between Model 1 and 2 can be seen from the development of the proxy i.e $V(\epsilon_{t-1}|y_{1:t-2})$ is replaced with $V(y_{t-1}|y_{1:t-2})$ in model 2.  
\subsubsection{Model 3}
In the previous approximation schemes, the information coming from $y_{t-1}$ is not used. \cite{Due97} suggests that $y_{t-1}$ should be included in the conditioning set of the proxy while assuming that $\mu_{t-1}$ and $\sigma_{t-1}^{2}$ are functions of $(y_{1:t-2},\xi_{t-2}')$. The following relation can thus be credited to him
\be
\begin{aligned}
\epsilon_{(D)t-1}&=y_{t-1}-\mu_{(D)t-1}\\
\mu_{(D)t-1}&=E[\mu_{t-1}(y_{1:t-2},\xi_{t-2}'\theta_{\mu})|y_{1:t-1}]=\sum_{m=1}^{M}\mu_{t-1}(y_{1:t-2},e_{m}'\theta_{\mu})q(\xi_{t-2}'=e_{m}'|y_{1:t-1})\\
\sigma_{(D)t-1}^{2}&= E[\sigma_{t-1}^{2}(y_{1:t-2},\xi_{t-2}'\theta_{\sigma})|y_{1:t-1}]			=\sum_{m=1}^{M}\sigma_{t-1}^{2}(y_{1:t-2},e_{m}'\theta_{\sigma})q(\xi_{t-2}'=e_{m}'|y_{1:t-1}).
\end{aligned}
\ee 
The probability $q(\xi_{t-1}'=e_{m}'|y_{1:t})$ is a one period ahead smoothed probability which can be computed as:
\be
\begin{aligned}
q(\xi_{t-1}'=e_{m}'|y_{1:t})&=\sum_{i=1}^{M}q(\xi_{t-1}'=e_{m}',\xi_{t}'=e_{i}'|y_{1:t})\\
		&=\sum_{i=1}^{M}q(\xi_{t-1}'=e_{m}'|\xi_{t}'=e_{i}',y_{1:t})q(\xi_{i}'=e_{i}'|y_{1:t})\\
		&=\sum_{i=1}^{M}q(\xi_{t-1}'=e_{m}'|\xi_{t}'=e_{i}',y_{1:t-1})q(\xi_{i}'=e_{i}'|y_{1:t})\\
		&=\sum_{i=1}^{M}\dfrac{q(\xi_{t-1}'=e_{m}',\xi_{t}'=e_{i}'|y_{1:t-1})q(\xi_{i}'=e_{i}'|y_{1:t})}{q(\xi_{t}'=e_{i}'|y_{1:t-1})}\\
		&=q(\xi_{t-1}'=e_{m}'|y_{1:t-1})\sum_{i=1}^{M}\dfrac{q(\xi_{t}'=e_{i}'|\xi_{t-1}'=e_{m}',y_{1:t-1})q(\xi_{i}'=e_{i}'|y_{1:t})}{q(\xi_{t}'=e_{i}'|y_{1:t-1})}
\end{aligned}
\ee
Within this framework we note that the conditioning is on $y_{1:t-1}$ while the functional form depends on  $(y_{1:t-2},\xi_{t-2}')$. We equally note that at every time step $t$ the value of $q(\xi_{t-2}'=e_{m}'|y_{1:t-1})$  for all $m$ is required for computation.
\subsubsection{Model 4}
The following approximation is similar to model 3. As opposed to model 3, we assume that $\mu_{t-1}$ and $\sigma_{t-1}^{2}$ are functions of $(y_{1:t-2},\xi_{t-1}')$. This modification leads to the following approximation  
 \cite{klas02} model.
\be
\begin{aligned}
\epsilon_{(SK)t-1}&=y_{t-1}-\mu_{(SK)t-1}\\
\mu_{(SK)t-1}&=E[\mu_{t-1}(y_{1:t-2},\xi_{t-1}'\theta_{\mu})|y_{1:t-1}]=\sum_{m=1}^{M}\mu_{t-1}(y_{1:t-2},e_{m}'\theta_{\mu})q(\xi_{t-1}'=e_{m}'|y_{1:t-1})\\
\sigma_{(SK)t-1}^{2}&= E[\sigma_{t-1}^{2}(y_{1:t-2},\xi_{t-1}'\theta_{\sigma})|y_{1:t-1}]					 =\sum_{m=1}^{M}\sigma_{t-1}^{2}(y_{1:t-2},e_{m}'\theta_{\sigma})q(\xi_{t-1}'=e_{m}'|y_{1:t-1}).
\end{aligned}
\ee 
In the next approximation, the current regime will be added to the conditioning set of this version of the auxiliary model. Hence, this approximation will be identified as the simplified version of \cite{klas02} model. In order to implement this approximation scheme the value of $q(\xi_{t-1}'=e_{m}'|y_{1:t-1})$ for all $m$ is required at each point in time $t$.
\subsubsection{Model 5}
In each of the approximations described above, information relating to the current regime is ignored in the conditioning set. On observing this, \cite{klas02} suggests the following approximation
\be
\begin{aligned}
\epsilon_{(K)t-1}&=y_{t-1}-\mu_{i,(K)t-1}\\
\mu_{i,(K)t-1}&=E[\mu_{t-1}(y_{1:t-2},\xi_{t-1}'\theta_{\mu})|y_{1:t-1},\xi_{t}'=e_{i}']\\
            &=\sum_{m=1}^{M}\mu_{t-1}(y_{1:t-2},e_{m}'\theta_{\mu})q(\xi_{t-1}'=e_{m}'|y_{1:t-1},\xi_{t}'=e_{i}')\\
\sigma_{i,(K)t-1}^{2}&= E[\sigma_{t-1}^{2}(y_{1:t-2},\xi_{t-1}'\theta_{\sigma})|y_{1:t-1}]\\
						 &=\sum_{m=1}^{M}\left(\mu_{t-1}(y_{1:t-2},e_{m}'\theta_{\mu})+\sigma_{t-1}^{2}(y_{1:t-2},e_{m}'\theta_{\sigma})\right)q(\xi_{t-1}'=e_{m}'|y_{1:t-1},\xi_{t}'=e_{i}')\\
						 &-\left(\sum_{m=1}^{M}\mu_{t-1}(y_{1:t-2},e_{m}'\theta_{\mu})q(\xi_{t-1}'=e_{m}'|y_{1:t-1},\xi_{t}'=e_{i}')\right)^{2},
\end{aligned}
\ee 
where
\be
\begin{aligned}
q(\xi_{t-1}'=e_{m}'|y_{1:t-1},\xi_{t}'=e_{i}')&=\dfrac{q(\xi_{t-1}'=e_{m}',\xi_{t}'=e_{i}'|y_{1:t-1})}{q(\xi_{t}'=e_{i}'|y_{1:t-1})}\\
		&=\dfrac{q(\xi_{t}'=e_{i}'|y_{1:t-1},\xi_{t-1}'=e_{m}')q(\xi_{t-1}'=e_{m}'|y_{1:t-1})}{q(\xi_{t}'=e_{i}'|y_{1:t-1})}
\end{aligned}
\ee
Note that this approximation requires the computation of $q(\xi_{t-1}'=e_{m}'|y_{1:t-1},\xi_{t}'=e_{i}')$ for all $m$ and $i$ at time $t$. 
\subsection{Sampling the $\theta$}
Sampling $\theta$ from the full conditional distribution will be done by separating the parameters of the transition matrix from the GARCH parameters. We assume that the parameters of the transition probabilities are independent of GARCH parameters.
\subsubsection{Sampling transition probability parameters}
The posterior distribution of $\theta_{\pi}$ is given by
\begin{equation}
\begin{aligned}
f(\theta_{\pi}|\xi_{1:T},\theta_{\mu},\theta_{\sigma},y_{1:T})&\propto f(\xi_{1:T},\theta_{\mu},\theta_{\sigma},y_{1:T}|\theta_{\pi})f(\theta_{\pi})\\
																					&\propto f(\xi_{1:T},y_{1:T}|\theta)f(\theta_{\pi})\\
																					&\propto f(\theta_{\pi})\prod_{t=2}^{T}\left(\prod_{i=1}^{M}(\pi_{i}\xi_{t-1})^{\xi_{it}}\right)\\
																					&=f(\theta_{\pi})\prod_{t=2}^{T}\left(\prod_{i=1}^{M}\left(\sum_{j=1}^{M}\pi_{ij}\xi_{jt-1}\right)^{\xi_{it}}\right)\\
																					&=f(\theta_{\pi})\prod_{j=1}^{M}\prod_{i=1}^{M}\pi_{ij}^{n_{ij}}
\end{aligned}\label{eq11}
\end{equation}
where $n_{ij}$ is the number of times $\xi_{it}=\xi_{jt-1}=1$ for $i,j=1,\dots,M$. It is easy to show that by substituting, as defined earlier, the conjugate Dirichlet prior for the transition probabilities, $\theta_{\pi}$, in (\ref{eq11}) we obtain
\begin{equation}
f(\theta_{\pi}|\xi_{1:T},\theta_{\mu},\theta_{\sigma},y_{1:T}) = \prod_{m=1}^{M}{\text{Dirichlet}}(n_{1m}+\eta_{1m},\dots,n_{Mm}+\eta_{Mm}).
\end{equation}
\subsubsection{Sampling GARCH parameters}
Given a prior density $f(\theta_{\mu},\theta_{\sigma})$, the posterior density of  $(\theta_{\mu},\theta_{\sigma})$ can be expressed as 
\begin{equation}
f(\theta_{\mu},\theta_{\sigma}|\xi_{1:T},\theta_{\pi},y_{1:T})\propto f(\theta_{\mu},\theta_{\sigma})\prod_{t=1}^{T}{\mathcal{N}}(\mu_{t}(y_{1:t-1},\xi_{t}'\theta_{\mu}),\sigma_{t}^{2}(y_{1:t-1},\xi_{t}'\theta_{\sigma}))\label{eq12}
\end{equation}
For this step of the Gibbs sampler, we apply adaptive Metropolis-Hastings (MH) sampling technique since the full conditional distribution is known to be non-standard. Details can be found, as required, in the appendix.
\section{Illustration with simulated data set} 
\begin{figure}[htp]
\centering
\includegraphics[width=1.0\textwidth, height = 70mm, clip]{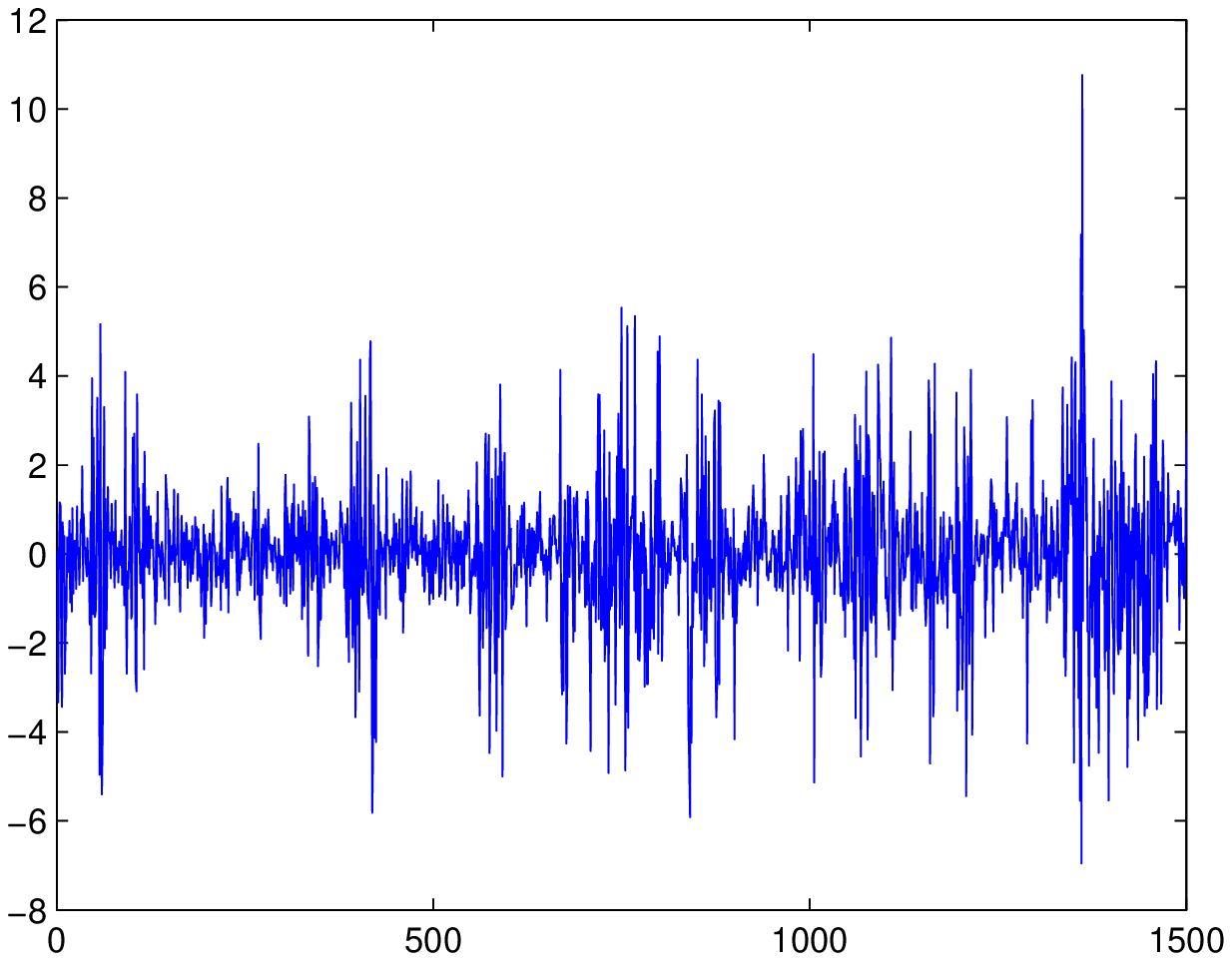}\\
\includegraphics[width=1.0\textwidth, height = 70mm,clip]{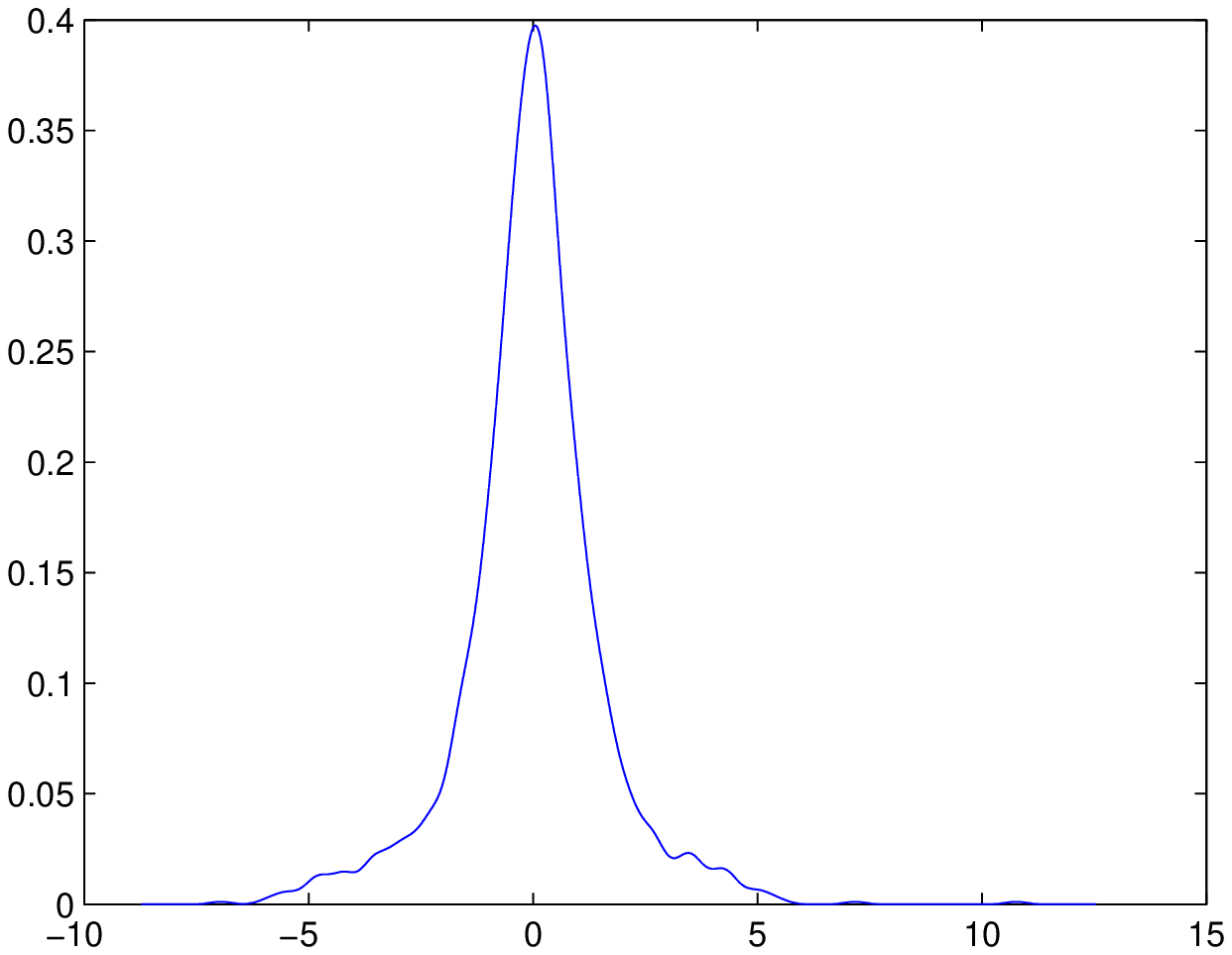}\\
\includegraphics[width=1.0\textwidth, height = 70mm,clip]{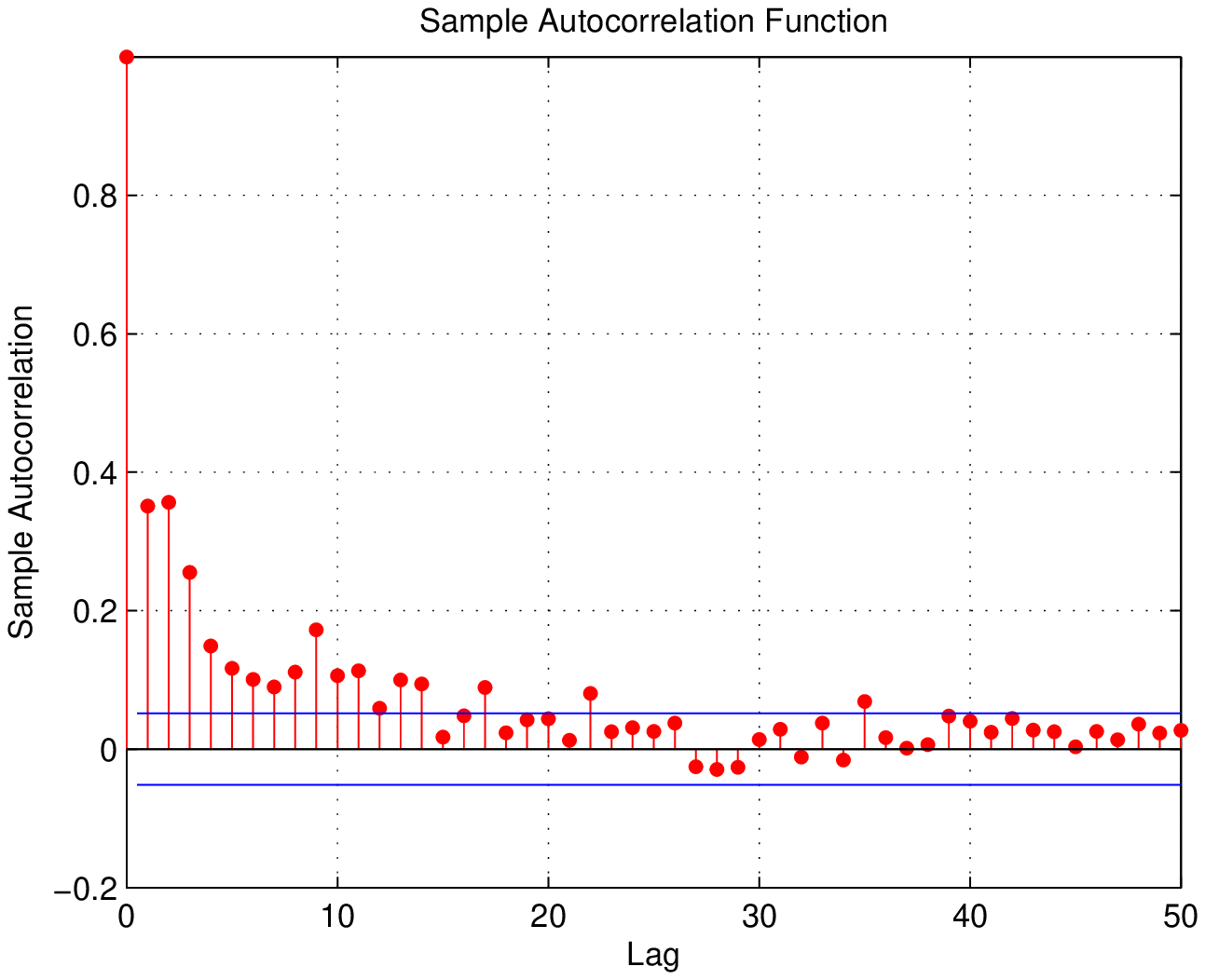}\\
\caption{Graphs for the simulated data for DGP defined in Table \ref{tab1}.}
\label{Fig1}
\end{figure}

We generated a time series of length 1500 from the data generating process corresponding to the model defined by equations (\ref{eq3}) and (\ref{eq4}) for two regimes ($M=2$), time invariant transition probabilities and constant parameter switching conditional mean. The parameter values for the simulation exercise are set at: $\mu =(\mu_{1},\mu_{2})=(0.06,-0.09)$, $\gamma =(\gamma_{1},\gamma_{2})=(0.30,2.00)$, $\alpha =(\alpha_{1},\alpha_{2})=(0.35,0.10)$, $\beta =(\beta_{1},\beta_{2})=(0.20,0.60)$, $\pi_{11}=0.98$, $\pi_{22}=0.96$. These parameter values corresponds to the choices made by \cite{Bau10} in a similar Monte Carlo exercise. A relatively higher and more persistent conditional variance as compared to the first GARCH equation is implied by the second GARCH equation. Also, the transition probabilities of remaining in each regime is close to one. A summary statistics of a typical series of length 1500 simulated from this DGP is reported in Table \ref{tab1} , and in Figure \ref{Fig1} we display, respectively, the time series, kernel density estimate and the autocorrelation function (ACF) of the square of the same series.  The mean of the series is close to zero and the excess kurtosis is estimated to be 3.57. 
\begin{table*}[!h]
	\caption{ Descriptive statistics for simulated data. }
	\centering
		\begin{tabular}{|cccccc|}
		\hline
		  Min.  &      max.   &    Mean     &     Std.  &   Skewness &    Kurtosis\\
  $-6.9540$ &    10.7600  &  $-0.0042$  &    1.5740 & $0.04120$  &     6.5659.\\
   	\hline
	\end{tabular}
	\label{tab1}
\end{table*}
For each hidden state sampling algorithm described in Section 3.1 and the auxiliary models presented in Section 3.2, we perform 10000 Gibbs iterations and compare estimates from these schemes with estimates obtained using the single-move sampling scheme for the hidden states. To carry out the MCMC exercise, we set the initial parameters of the algorithm to the maximum likelihood estimates of one of the MS-GARCH approximations described in Section 3.1 and randomly generated initial state trajectory.  The hyperparameters of the prior distributions of the transition probabilities $\nu_{ij}$ for $i,j=1,2$ are set to 1 while the support for other parameters are given in the table reporting their parameter estimates. The case of two trials, $(K=2)$, is considered within the different multi-point sampling strategies discussed earlier. Table from \ref{tab:Model1} to  \ref{tab:Model5} highlights the posterior means and standard deviation of the parameters and the transition probabilities of the MS-GARCH under each of the auxiliary models used in constructing proposals for the hidden states. Column 4 of each of these tables reports the parameter estimates and transition probabilities obtained by using the single move technique for sampling the state variables within the Gibbs algorithm while in columns 5 to 7 we present the result obtained using the different multi-move multipoint sampling techniques within the Gibbs algorithm.
\begin{table*}[!h]
	\caption{ Estimated parameter value and posterior statistics using Model 1.}
	\centering
	\begin{tabular}{|c |c | c | c | c| c| c|}
\hline
             	     &            &                &                 &  \multicolumn{3}{c|}{Multi move}\\
\cline{5-7}			                 
		        	     & DGP Values & Prior supports &  Single Move 	 &    MTM  	 &	MTMIS 	  &   	MCTM  \\
		\hline
    $\pi_{11}$     &   0.980    &  (0.00 1.00)   &    0.968        &   0.972   &  0.974     &    0.977 \\
                   &            &                &   (0.014)       &  (0.005)  &  (0.006)   &  (0.005)\\
    $\pi_{22}$     &   0.960    &  (0.00 1.00)   &    0.995        &   0.952   &  0.955     &    0.957 \\
                   &            &                &    (0.002)      &   (0.011) &  (0.011)   &  (0.009) \\
    $\mu_{1}$      &   0.060    &  (0.02 0.15)   &    0.099        &   0.045   &  0.049     &    0.046 \\
                   &            &                &    (0.031)      &   (0.017) &  (0.019)   &   (0.0173)\\
    $\mu_{2}$      & $-0.090$   &($-0.35$ 0.18)  &  $-0.013$       & $-0.109$  &$-0.107$    &  $-0.110$\\
                   &            &                &	 (0.035)       &   (0.106) &  (0.108)   &   (0.107)\\
    $\gamma_{1}$   &   0.300    &  (0.15 0.45)   &    0.290        &   0.345   &  0.365     &    0.350 \\
                   &            &                &    (0.053)      &   (0.046) &  (0.046)   &   (0.047)\\
    $\gamma_{2}$   &   2.000    &  (0.50 4.00)   &    0.508        &   1.682   &  2.042     &    2.533 \\
                   &            &                & 		(0.010)      &   (0.432) &  (0.599)   &    (0.650)\\
    $\alpha_{1}$   &   0.350    &  (0.10 0.50)   &    0.227        &   0.141   &  0.181     &    0.180 \\
                   &            &                &    (0.099)      &   (0.037) &  (0.049)   &    (0.044)\\
    $\alpha_{2}$   &   0.100    &  (0.02 0.35)   &    0.331        &   0.042   &  0.047     &    0.047\\
                   &            &                &    (0.016)      &   (0.019) &  (0.023)   &    (0.024)\\
    $\beta_{1}$    &   0.200    &  (0.05 0.40)   &    0.190        &   0.248   &  0.196     &    0.227 \\ 
                   &            &                &    (0.097)      &   (0.082) &  (0.076)   &    (0.079)\\
    $\beta_{2}$    &   0.600    &  (0.35 0.85)   &    0.510        &   0.683   &  0.612     &    0.534 \\
                   &            &                &    (0.019)      &   (0.084) &  (0.109)   &    (0.111)\\
    \hline
    \end{tabular}    
	\label{tab:Model1}
\end{table*}
\begin{table*}[!h]
	\caption{ Estimated parameter value and posterior statistics using Model 2.}
	\centering
	\begin{tabular}{|c |c | c | c | c| c| c|}
\hline
             	     &            &                &                 &  \multicolumn{3}{c|}{Multi move}\\
\cline{5-7}			                 
		        	     & DGP Values & Prior supports &  Single Move 	 &    MTM  	 &	MTMIS 	  &   	MCTM  \\
		\hline
    $\pi_{11}$     &   0.980    &  (0.00 1.00)   &    0.968        &   0.973   &  0.9753    &   0.9771 \\
                   &            &                &    (0.014)      &  (0.006)  &  (0.006)   &   (0.006)\\                   
    $\pi_{22}$     &   0.960    &  (0.00 1.00)   &    0.995        &   0.952   &  0.952     &   0.957 \\
                   &            &                &    (0.002)      &  (0.011)  &  (0.011)   &   (0.010)\\
    $\mu_{1}$      &   0.060    &  (0.02 0.15)   &    0.099        &   0.045   &  0.047     &   0.048 \\
                   &            &                &    (0.031)      &  (0.017)  &  (0.018)   &   (0.018)\\
    $\mu_{2}$      & $-0.090$   &($-0.35$ 0.18)  &  $-0.013$       & $-0.108$  &$-0.111$    &   $-0.120$\\
                   &            &                &    (0.035)      &  (0.107)  &  (0.111)   &   (0.109)\\
    $\gamma_{1}$   &   0.300    &  (0.15 0.45)   &    0.290        &   0.344   &  0.328     &   0.347 \\
                   &            &                &    (0.052)      &  (0.046)  &  (0.052)   &   (0.052)\\
    $\gamma_{2}$   &   2.000    &  (0.50 4.00)   &    0.508        &   1.701   &  1.923     &   1.968 \\
                   &            &                &    (0.009)      &   (0.442) &  (0.626)   &   (0.673)\\
    $\alpha_{1}$   &   0.350    &  (0.10 0.50)   &    0.228        &   0.142   &  0.181     &   0.186 \\
                   &            &                &    (0.098)      &  (0.039)  &  (0.042)   &   (0.044)\\
    $\alpha_{2}$   &   0.100    &  (0.02 0.35)   &    0.331        &   0.042   &  0.043     &   0.044\\
                   &            &                &    (0.016)      &  (0.019)  & (0.021)    &  (0.022)\\
    $\beta_{1}$    &   0.200    &  (0.05 0.40)   &    0.190        &   0.250   &  0.275     &   0.237 \\ 
                   &            &                &    (0.096)      &   (0.079) &  (0.084)   &  (0.086)\\
    $\beta_{2}$    &   0.600    &  (0.35 0.85)   &    0.511        &   0.681   &  0.645     &   0.631 \\
                   &            &                &    (0.019)      &  (0.085)  &  (0.117)   &  (0.1216)\\
    \hline
    \end{tabular}    
	\label{tab:Model2}
\end{table*}
\begin{table*}[!h]
	\caption{ Estimated parameter value and posterior statistics using Model 3.}
	\centering
	\begin{tabular}{|c |c | c | c | c| c| c|}
\hline
             	     &            &                &                 &  \multicolumn{3}{c|}{Multi move}\\
\cline{5-7}			                 
		        	     & DGP Values & Prior supports &  Single Move 	 &    MTM  	 &	MTMIS 	  &   	MCTM  \\
		\hline
    $\pi_{11}$     &   0.980    &  (0.00 1.00)   &    0.968        &   0.975   &  0.976     &   0.977\\
    							 &            &                &    (0.014)      &   (0.005) &  (0.006)   &   (0.006)\\
    $\pi_{22}$     &   0.960    &  (0.00 1.00)   &    0.995        &   0.956   &  0.956     &   0.956 \\
    							 &            &                &    (0.002)      &   (0.009) &  (0.011)   &   (0.011)\\
    $\mu_{1}$      &   0.060    &  (0.02 0.15)   &    0.099        &   0.050   &  0.050     &   0.049 \\
    							 &            &                &    (0.031)      &   (0.018) &  (0.019)   &   (0.018)\\
    $\mu_{2}$      & $-0.090$   &($-0.35$ 0.18)  &  $-0.013$       & $-0.128$ &$-0.122$     &   $-0.116$\\
    							 &            &                &    (0.034)      &   (0.104) &  (0.106)   &   (0.108)\\
    $\gamma_{1}$   &   0.300    &  (0.15 0.45)   &    0.290        &   0.382   &  0.371     &   0.354 \\
    							 &            &                &    (0.052)      &   (0.043) &  (0.046)   &   (0.051)\\
    $\gamma_{2}$   &   2.000    &  (0.50 4.00)   &    0.508        &   2.107   &  2.059     &   2.448\\
    							 &            &                &    (0.009)      &   (0.641) &  (0.648)   &   (0.712)\\
    $\alpha_{1}$   &   0.350    &  (0.10 0.50)   &    0.227        &   0.168   &  0.174     &   0.167 \\
    							 &            &                &    (0.098)      &   (0.042) &  (0.047)   &   (0.047)\\
    $\alpha_{2}$   &   0.100    &  (0.02 0.35)   &    0.331        &   0.046   &  0.046     &   0.048\\
    							 &            &                &    (0.016)      &   (0.023) &  (0.022)   &   (0.025)\\
    $\beta_{1}$    &   0.200    &  (0.05 0.40)   &    0.190        &   0.173   &  0.199     &   0.237 \\ 
    							 &            &                &   (0.096)       &  (0.076)  &  (0.081)   &   (0.089)\\
    $\beta_{2}$    &   0.600    &  (0.35 0.85)   &    0.510        &   0.603   &  0.613     &   0.547 \\
    							 &            &                &    (0.019)      &   (0.114) &  (0.117)   &   (0.119)\\
    \hline
    \end{tabular} 
	\label{tab:Model3}
\end{table*}
\begin{table*}[!h]
	\caption{ Estimated parameter value and posterior statistics using Model 4.}
	\centering
	\begin{tabular}{|c |c | c | c | c| c| c|}
\hline
             	     &            &                &                 &  \multicolumn{3}{c|}{Multi move}\\
\cline{5-7}			                 
		        	     & DGP Values & Prior supports &  Single Move 	 &    MTM  	 &	MTMIS 	  &   	MCTM  \\
		\hline
    $\pi_{11}$     &   0.980    &  (0.00 1.00)   &    0.968        &   0.978   &  0.977     &   0.977 \\
                   &            &                &    (0.014)      &  (0.005)  & (0.006)    &  (0.005)\\
    $\pi_{22}$     &   0.960    &  (0.00 1.00)   &    0.995        &   0.959   &  0.958     &   0.957 \\
                    &            &               &    (0.002)      &  (0.010)  & (0.010)    &  (0.011)\\
    $\mu_{1}$      &   0.060    &  (0.02 0.15)   &    0.099        &   0.049   &  0.048     &   0.050 \\
                    &            &               &    (0.031)      &  (0.019)  & (0.018)    &  (0.019)\\
    $\mu_{2}$      & $-0.090$   &($-0.35$ 0.18)  &  $-0.013$       & $-0.121$  &$-0.117$    & $-0.134$\\
                    &            &               &    (0.034)      &  (0.109)  &  (0.108)   &  (0.108)\\
    $\gamma_{1}$   &   0.300    &  (0.15 0.45)   &    0.290        &   0.362   &  0.366     &   0.370 \\
                    &            &               &    (0.052)      &   (0.045) &  (0.046)   &  (0.0469)\\
    $\gamma_{2}$   &   2.000    &  (0.50 4.00)   &    0.508        &   2.519   &  1.931     &   2.173 \\
                    &            &               &    (0.009)      &   (0.683) &  (0.648)   &  (0.665)\\
    $\alpha_{1}$   &   0.350    &  (0.10 0.50)   &    0.227        &   0.170   &  0.179     &   0.172 \\
                    &            &               &    (0.098)      &   (0.041) &  (0.050)   &  (0.044)\\
    $\alpha_{2}$   &   0.100    &  (0.02 0.35)   &    0.331        &   0.046   &  0.046     &   0.046\\
                    &            &               &    (0.016)      &  (0.023)  &  (0.022)   &   (0.023)\\
    $\beta_{1}$    &   0.200    &  (0.05 0.40)   &    0.190        &   0.230   &  0.204     &   0.205 \\ 
                    &            &               &    (0.096)      &  (0.082)  &  (0.077)   &   (0.082)\\
    $\beta_{2}$    &   0.600    &  (0.35 0.85)   &    0.510        &   0.539   &  0.633     &   0.594 \\
                    &            &               &    (0.019)      &   (0.113) &  (0.116)   &   0.1157\\
    \hline
    \end{tabular} 
	\label{tab:Model4}
\end{table*}
\begin{table*}[!h]
	\caption{ Estimated parameter value and posterior statistics using Model 5.}
	\centering
	\begin{tabular}{|c |c | c | c | c| c| c|}
\hline
             	     &            &                &                 &  \multicolumn{3}{c|}{Multi move}\\
\cline{5-7}			                 
		        	     & DGP Values & Prior supports &  Single Move 	 &    MTM  	 &	MTMIS 	  &   	MCTM  \\
		\hline
    $\pi_{11}$     &   0.980   &  (0.00 1.00)    &    0.968        &   0.974   &  0.976     &   0.976 \\
                    &            &               &   (0.015)       &   (0.006) & (0.006)    &  (0.006)\\ 
    $\pi_{22}$     &   0.960   &  (0.00 1.00)    &    0.995        &   0.954   &  0.957     &   0.957 \\
                    &            &               &   (0.002)       &  (0.012)  & (0.011)    &  (0.011) \\
    $\mu_{1}$      &   0.060   &  (0.02 0.15)    &    0.099        &   0.050   &  0.049     &   0.050\\
                    &           &                &   (0.031)       &  (0.019)  &  (0.018)   &   (0.019)\\ 
    $\mu_{2}$      & $-0.090$  &($-0.35$ 0.18)   &  $-0.013$       & $-0.127$  & $-0.124$   &  $-0.123$\\
                    &            &               &  (0.035)        &   (0.107) & (0.108)    &  (0.105)\\
    $\gamma_{1}$   &   0.300   &  (0.15 0.45)    &    0.290        &   0.368   &  0.373     &   0.378 \\
                    &            &               &  (0.053)        &   (0.045) &  (0.046)   & (0.045)\\ 
    $\gamma_{2}$   &   2.000   &  (0.50 4.00)    &    0.508        &   1.869   &  1.864     &   2.069 \\
                    &            &               &   (0.010)       &   (0.694) &  (0.679)   &  (0.629)\\  
    $\alpha_{1}$   &   0.350   &  (0.10 0.50)    &    0.227        &   0.172   &  0.171     &   0.177 \\
                    &            &               &  (0.098)        &   (0.044) &  (0.044)   &   (0.046)\\
    $\alpha_{2}$   &   0.100   &  (0.02 0.35)    &    0.331        &   0.045   &  0.045     &   0.047\\
                    &            &               &  (0.016)        &   (0.022) &  (0.022)   &   (0.024)\\
    $\beta_{1}$    &   0.200   &  (0.05 0.40)    &    0.190        &   0.200   &  0.194     &   0.183 \\ 
                    &            &               &  (0.096)        &  (0.079)  &  (0.079)   &   (0.079)\\
    $\beta_{2}$    &   0.600   &  (0.35 0.85)    &    0.510        &   0.648   &  0.648     &   0.608 \\
                    &            &               &  (0.019)        &  (0.126)  &  (0.123)   &   (0.116)\\
    \hline
    \end{tabular}
	\label{tab:Model5}
\end{table*}
\begin{figure}[htb]
	\centering
	\includegraphics[width=1.0\textwidth, height = 100mm]{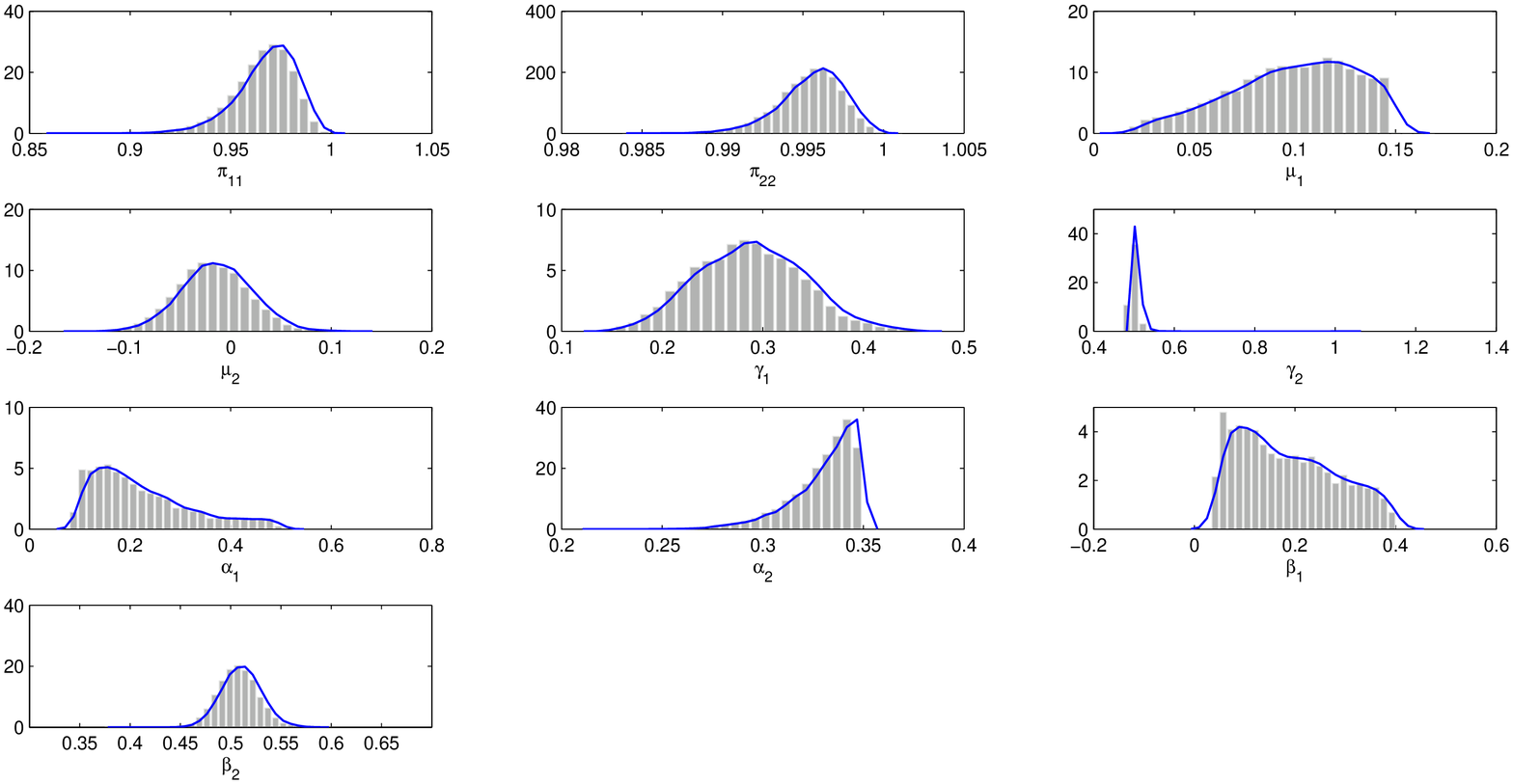}
	\caption{Posterior densities of the MS-GARCH parameters using single-move Scheme}
	\label{fig:singlemove}
\end{figure}
\begin{figure}[htb]	 
	\includegraphics[width=1.0\textwidth, height = 100mm]{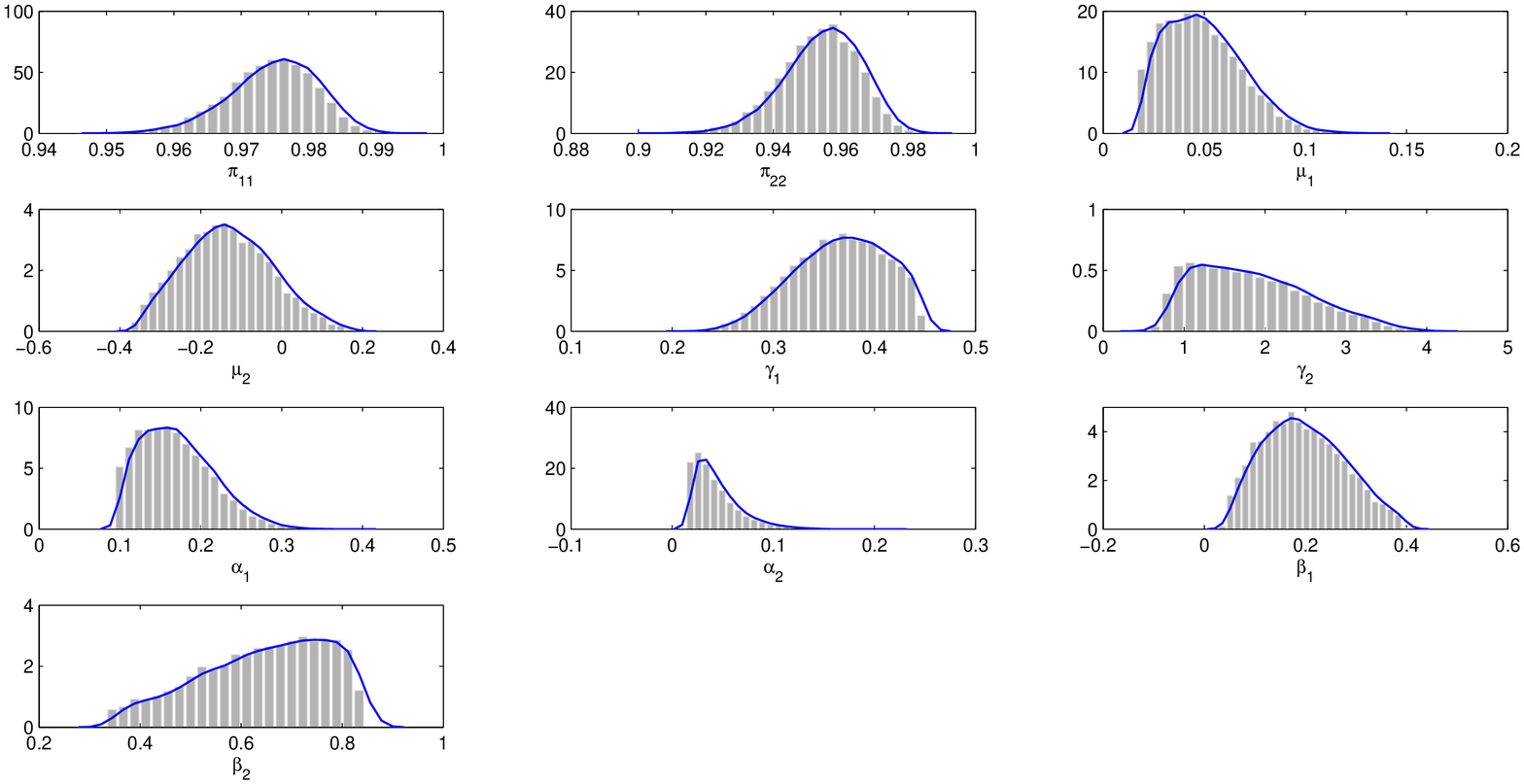}
	\caption{Posterior densities of the MS-GARCH parameters using MTM with model 5}
	\label{fig:MTM05}
\end{figure}
\begin{figure}[htb]	 
	\includegraphics[width=1.0\textwidth, height = 100mm]{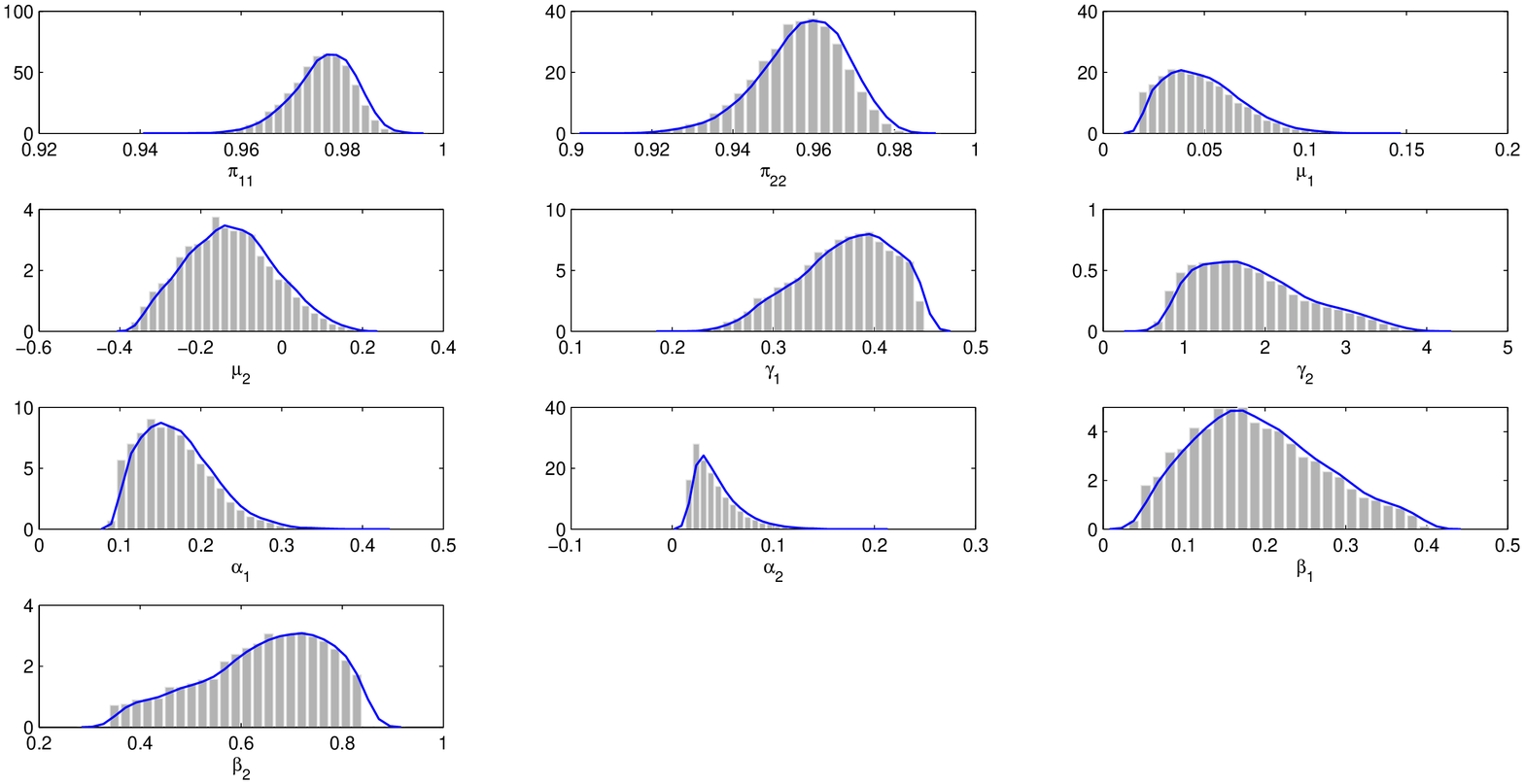}
	\caption{Posterior densities of the MS-GARCH parameters using MTMIS with model 5}
	\label{fig:MTMIS05}
\end{figure}
\begin{figure}[htb]	 
	\includegraphics[width=1.0\textwidth, height = 100mm]{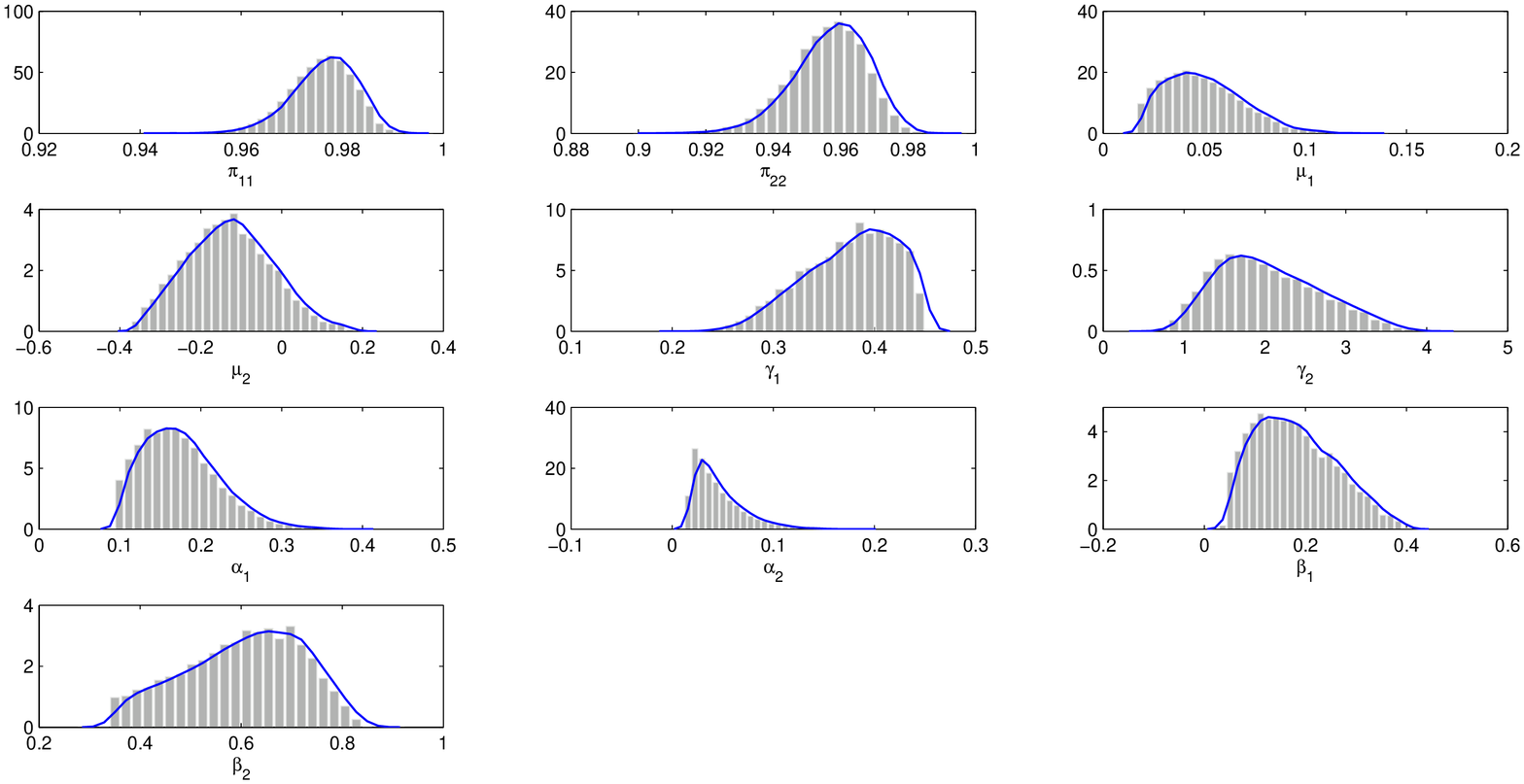}
	\caption{Posterior densities of the MS-GARCH parameters using MCTM with model 5}
	\label{fig:MTCM05}
\end{figure}
With the exception of a few, the posterior means under the multi-move multi-point sampling schemes relative to the single-move technique have more values within one posterior standard deviation away from the DGP values. In Figure from \ref{fig:singlemove} to \ref{fig:MTCM05} we report the posterior densities of the parameters using single-move, MTM, MTMIS, and MTCM sampling strategies respectively. The multi-move sampler are constructed using model 5. The shapes of the posterior densities are unimodal, thus ruling out label switching problem. We also examine the performance of our multi-move multipoint algorithms relative to the single-move strategy by computing the percentage of correctly specified regimes. To do this, we first calculate the average of the Gibbs output on the state variables and then assign mean states greater than one-half to regime 2 (and regime 1 otherwise). We find out that the single-move technique is able to classify $43\%$ of the data correctly while the multi-move multipoint samplers classified between $93\%$
and $96\%$ of the data correctly. The acceptance rate of the the multi-move multipoint proposals varies between $1\%$ and $20\%$ with the highest arising from multipoint sampling schemes proposal distribution constructed using model 5. We compute the mean squared error (MSE) of the posterior means of parameter relative to the true parameter to further quantify our estimators i.e.
\begin{equation}
MSE = \dfrac{1}{n}\sum_{i=1}^{n}(\hat{\theta}_{i} - \theta_{i})^{2}
\end{equation}
where $n$ is the number of parameters, $\hat{\theta}_{i}$ is the parameter estimate of the $i$-th element, $\theta_{i}$, of the DGP parameter set. The result of this exercise is reported in Table {\ref{tab7}. From Table {\ref{tab7}}, the low MSE of our multipoint sampling schemes further confirms their superiority over the single-move procedures. 
\begin{table*}[!h]
	\caption{ Mean Squared Error (MSE). }
	\centering
		\begin{tabular}{|ccccc|}
		\hline
		          &    Single move   &    MTM     &     MTMIS    &   MCTM\\
		\hline
     Model 1  &      0.2310      &  0.0160    &    0.0038    &  0.0324\\
     Model 2  &      0.2310      &  0.0147    &    0.0047    &  0.0036\\
     Model 3  &      0.2310      &  0.0056    &    0.0044    &  0.0245\\
     Model 4  &      0.2310      &  0.0315    &    0.0043    &  0.0071\\
     Model 5  &      0.2310      &  0.0060    &    0.0062    &  0.0045\\
   	\hline
	\end{tabular}
	\label{tab7}
\end{table*}
The inefficiency of the various multi-move multiple-try Metropolis samplers relative to the single-move sampler are further assessed by examining how much the variance of the parameters are increased due to autocorrelation coming from the sampler. Let $z^{(1)},\dots,z^{(G)}$ denote a sample from the posterior distribution of a random variable $Z$. Then inefficiency factor $(IF)$ is evaluated by 
\begin{equation}
IF = 1 + 2\sum_{l=1}^{L}w_{l}\rho_{l}
\end{equation}
where $\rho_{l}$, $l=1,2,\dots$ is the autocorrelation function of $z^{(1)},\dots,z^{(G)}$ at lag $l$ and $w_{l}$ is the associated weight. If the samples are independent then $IF = 1$. If $A$ and $B$ are two competing algorithm with inefficient factor $IF_{A}$ and $IF_{B}$ respectively then we define the relative inefficiency $(RI)$ as:
\begin{equation}
RI = \dfrac{Time_{A}}{Time_{B}}\times\dfrac{IF_{A}}{IF_{B}}
\end{equation}
where $Time_{A}$ and $Time_{B}$ corresponds to the computing times of each algorithm. $RI$ measures the factor by which the run-time of algorithm $A$ must be increased to achieve algorithm $B$'s precision; values greater than one suggests that algorithm $B$ is more efficient. We provide in Table from to \ref{tab:RImodel5} the $RI$ for various multi-move multipoint algorithms relative to the single-move sampling strategy. The number of lags over which we calculate the $RI$ is fixed at $L=500$. From these tables our multi-move multipoint algorithms are more efficient than the single-move sampling technique for the state variable. This is despite the low acceptance rate of the of the multipoint proposals. Finally we shall notice that, as discussed in \cite{Craiu:2007rr}, a larger number of proposals is required to observe an appreciable difference in the efficiency of the MCTM over the standard MTM. 
\begin{table*}[!h]
	\caption{ Relative inefficiency factor using Model 1}
	\centering
		\begin{tabular}{|cccc|}
		\hline
		                              &    MTM      &     MTMIS    &   MCTM\\
		\hline
		$\max_{t=1:T}(\sigma^{2}_{t})$&   68.16   &   95.79    &   92.48\\
		  $\pi_{11}$                  &   64.31   &   60.37    &  139.11\\
      $\pi_{22}$                  &   53.93   &   65.52    &  115.91\\
      $\mu_{1}$                   &   119.42  &   105.59   &  153.58\\
      $\mu_{2}$                   &   78.08   &    62.13   &  107.04\\
      $\gamma_{1}$                &   45.96   &    77.43   &   66.18\\
      $\gamma_{2}$                &   14.69   &    17.57   &  15.29\\
      $\alpha_{1}$                &   77.54   &   136.39   &  206.11\\
      $\alpha_{2}$                &   42.54   &   64.04    &  71.15\\
      $\beta_{1}$                 &   44.76   &   89.79    &  76.29\\
      $\beta_{2}$                 &   26.05   &   32.98    &  29.96\\     
   	\hline
	\end{tabular}
	\label{tab:RImodel1}
\end{table*}
\begin{table*}[!h]
	\caption{ Relative inefficiency factor using Model 2}
	\centering
		\begin{tabular}{|cccc|}
		\hline
		                              &    MTM      &     MTMIS    &   MCTM\\
		\hline
		$\max_{t=1:T}(\sigma^{2}_{t})$&   72.08     &     93.97    &   95.35\\
		  $\pi_{11}$                  &   54.26     &     71.36    &   82.63\\
      $\pi_{22}$                  &   53.43     &     60.85    &   86.16\\
      $\mu_{1}$                   &   125.27    &     124.69   &   156.10\\
      $\mu_{2}$                   &   81.05     &     78.37    &   66.96\\
      $\gamma_{1}$                &   50.08     &     53.53    &   55.99\\
      $\gamma_{2}$                &   15.11     &     16.20    &   14.21\\
      $\alpha_{1}$                &   76.74     &     238.36   &   202.02\\
      $\alpha_{2}$                &   45.30     &     58.35    &   60.34\\
      $\beta_{1}$                 &   49.03     &     62.00    &   63.08\\
      $\beta_{2}$                 &   26.94     &     28.97    &   26.60  \\ 
   	\hline
	\end{tabular}
	\label{tab:RImodel2}
\end{table*}
\begin{table*}[!h]
	\caption{ Relative inefficiency factor using Model 3}
	\centering
		\begin{tabular}{|cccc|}
		\hline
		                              &    MTM      &     MTMIS    &   MCTM\\
		\hline
		$\max_{t=1:T}(\sigma^{2}_{t})$&   66.64     &     94.80    &  90.29\\
		  $\pi_{11}$                  &   55.04     &     51.68    &  58.42\\
      $\pi_{22}$                  &   63.59     &     62.76    &  49.31\\
      $\mu_{1}$                   &   96.03     &    107.90    & 147.03\\
      $\mu_{2}$                   &   50.53     &     71.94    &  84.67\\
      $\gamma_{1}$                &   49.08     &     72.63    &  55.65\\
      $\gamma_{2}$                &   10.64     &     15.14    &  13.64\\
      $\alpha_{1}$                &   129.17    &    142.76    & 114.75\\
      $\alpha_{2}$                &   39.85     &     60.02    &  61.12\\
      $\beta_{1}$                 &   50.69     &     75.29    &  59.40\\
      $\beta_{2}$                 &   19.97     &     28.55    &  26.43\\    
   	\hline
	\end{tabular}
	\label{tab:RImodel3}
\end{table*}
\begin{table*}[!h]
	\caption{ Relative inefficiency factor using Model 4}
	\centering
		\begin{tabular}{|cccc|}
		\hline
		                              &    MTM      &     MTMIS    &   MCTM\\
		\hline
		$\max_{t=1:T}(\sigma^{2}_{t})$&   74.01     &    96.79     &  94.01\\
		  $\pi_{11}$                  &   44.37     &    62.01     &  77.53\\
      $\pi_{22}$                  &   68.24     &    76.50     &  59.64\\
      $\mu_{1}$                   &   97.07     &    156.67    &  142.73\\
      $\mu_{2}$                   &   60.36     &    71.65     &  50.81\\
      $\gamma_{1}$                &   58.35     &    75.87     &  65.73\\
      $\gamma_{2}$                &   11.15     &    15.45     &  15.35\\
      $\alpha_{1}$                &   174.85    &    129.64    &  180.54\\
      $\alpha_{2}$                &   50.28     &    59.96     &  63.24\\
      $\beta_{1}$                 &   53.23     &    83.88     &  68.51\\
      $\beta_{2}$                 &   22.25     &    28.95     &  29.81  \\   
   	\hline
	\end{tabular}
	\label{tab:RImodel4}
\end{table*}
\begin{table*}[!h]
	\caption{ Relative inefficiency factor using Model 5}
	\centering
		\begin{tabular}{|cccc|}
		\hline
		                              &    MTM      &     MTMIS    &   MCTM\\
		\hline
		$\max_{t=1:T}(\sigma^{2}_{t})$&   69.05     &   92.88      &  114.51\\
		  $\pi_{11}$                  &   41.02     &   71.41      &   64.78\\
      $\pi_{22}$                  &   47.10     &   73.47      &   69.97\\
      $\mu_{1}$                   &   96.93     &  135.98      &  157.25\\
      $\mu_{2}$                   &   46.60     &   67.22      &   81.80\\
      $\gamma_{1}$                &   55.87     &   75.55      &   80.21\\
      $\gamma_{2}$                &   9.39      &   14.55      &   16.68\\
      $\alpha_{1}$                &   125.95    &  185.58      &  179.61\\
      $\alpha_{2}$                &   41.49     &   57.63      &   56.37\\
      $\beta_{1}$                 &   57.95     &   83.33      &   85.43\\
      $\beta_{2}$                 &   17.35     &   26.76      &   30.32 \\    
   	\hline
	\end{tabular}
	\label{tab:RImodel5}
\end{table*}
\clearpage
\section{Conclusion}
In this paper we deal with the challenging issue of efficient sampling algorithm for Bayesian inference on Markov-switching GARCH models. We provide some new algorithms based on the combination of multi-move and multi-points strategies.

More specifically, we apply the multiple-try sampler of \cite{Craiu:2007rr} combined with multi-move Gibbs sampler to Markov-switching GARCH models. For generating correlated proposal, we introduce antithetic Forward Filtering Backward Sampling (FFBS) algorithm for MS-GARCH based on the permuted displacement method of \cite{crameng2}. Our algorithms also extend to Markov-switching state space models the algorithms of \cite{So06} for continuous state space models.

From the results of our computational exercise, we observed a substantial gain in the efficiency of our Gibbs samplers over the usual single-move sampling algorithm for estimating the parameters of the MS-GARCH model. We also observed low acceptance rate $(1\% - 20\%)$ for the multipoint proposals. Despite the low acceptance rate for the multipoint proposals, we still have good results considering the length of the time series (1500) used. We expect that using the blocking scheme (as in  \cite{So06}) the efficiency and the acceptance rate of can our sampling procedure may increase. The issues of the choice of block length and of the application of the inference procedure to real data could be a matter of future research.
\clearpage
\section*{Appendix}
\subsection*{Constructing proposal distribution for $\theta_{\mu},\theta_{\sigma}$}
 Sample $\theta_{\mu}^{(r)},\theta_{\sigma}^{(r)}$ from $f(\theta_{\mu},\theta_{\sigma}|\xi_{1:T}^{(r)},\pi^{(r)},y_{1:T})$. Given a prior density $f(\theta_{\mu},\theta_{\sigma})$, the posterior density of  $\theta_{\theta_{\mu},\theta_{\sigma}}$ can be expressed as follows
\begin{equation}
f(\theta_{\mu},\theta_{\sigma}|\xi_{1:T}^{(r)},\pi,y_{1:T})\propto f(\theta_{\mu},\theta_{\sigma})\prod_{t=1}^{T}{\mathcal{N}}(y_{t};{\xi_{t}^{(r)}}'{\mu},\sigma_{t}^{2})
\end{equation}
where,
\be
\sigma_{t}^{2}={\xi_{t}^{(r)}}'{\gamma}+({\xi_{t}^{(r)}}'{\alpha})(y_{t-1}-{\xi_{t-1}^{(r)}}'\mu)^{2}+({\xi_{t}^{(r)}}'{\beta})\sigma_{t-1}^{2}.
\ee

In order to generate $\theta_{\mu},\theta_{\sigma}$ from the joint distribution we apply a further blocking of the Gibbs sampler. First, in the spirits of \cite{Fru06} we consider the full conditional distributions of the regime-specific parameters, and secondly, we split the regime-dependent parameters in two subvectors, the parameter of the observation equation and the parameters of the volatility process.
As regards the parameters of the return process equation, 
\be
f(\mu_{k}|\xi_{1:T}^{(r)},\mu_{-k}^{(r-1)},\gamma^{(r-1)},\beta^{(r-1)},\alpha^{(r-1)},y_{1:T})\propto \prod_{t\in\mathcal{T}_{k}}{\mathcal{N}}(y_{t};\mu_{k},\sigma_{t}^{2})\prod_{t\in\mathcal{T}_{k}^{-}}{\mathcal{N}}(y_{t};{\xi_{t}^{(r)}}'\mu,\sigma_{t}^{2})
\ee
where $\mu_{-k}=(\mu_{1},\ldots,\mu_{k-1},\mu_{k+1},\ldots,\mu_{M})'$, $\mathcal{T}_{k}=\{t=1,\ldots,T|\xi_{k,t}^{(r)}=1\}$, $\mathcal{T}_{k}^{-}=\{t=1,\ldots,T|\xi_{k,t}^{(r)}=0,\xi_{k,t-1}^{(r)}=1\}$.
It is not possible to simulate exactly from the full conditional distribution of $\mu_{k}$, $k=1,\ldots,M$ given the other parameters and the allocation variables, thus we apply a MTM step with independent normal proposal distribution. Focusing on the first term of the full conditional
\be
\prod_{t\in \mathcal{T}_{k}}\frac{1}{\sqrt{2\pi\sigma_{t}^{2}}}\exp\left\{-\frac{1}{2}\left(\mu_{k}^2 \sum_{t\in \mathcal{T}_{k}}\sigma_{t}^{-2}-2\mu_{k}\sum_{t\in \mathcal{T}_{k}}y_{t}\sigma_{t}^{-2}+\sum_{t\in \mathcal{T}_{k}}y_{t}^{2}\sigma_{t}^{-2}\right)\right\}
\ee
and if an approximation $\sigma^{*2}_{t}$ of $\sigma_{t}^{2}$  is available, then it is possible to approximate this part of the full conditional with a normal distribution with mean and variance
$$
m_{k}=s^{2}_{k}\left(\sum_{t\in\mathcal{T}_{k}}y_{t}/\sigma^{*2}_{t}\right),\quad s^{2}_{k}=\left(\sum_{t\in\mathcal{T}_{k}}1/\sigma^{*2}_{t}\right)^{-1}
$$
respectively, where $$\sigma_{t}^{*2}=({\xi_{t}^{(r)}}'\gamma^{(r-1)})+({\xi_{t}^{(r)}}'\alpha^{(r-1)})(y_{t-1}-{\xi_{t-1}^{(r)}}'\mu^{*})^{2}+({\xi_{t}^{(r)}}'\beta^{(r-1)})\sigma_{t-1}^{*2}$$ 
with $\mu^{*}=(\mu_{1}^{*},\ldots,\mu_{M}^{*})$, $\mu_{j}^{*}=T_{j}^{-1}\sum_{t\in\mathcal{T}_{j}}y_{t}$ and $T_{j}=\sum_{t\in\mathcal{T}_{j}}\xi_{j,t}$. This normal can be used as proposal in the MH step.

As regards the parameters of the volatility process the full conditional is
\begin{equation}
f(\gamma_{k},\beta_{k},\alpha_{k}|\xi_{1:T}^{(r)},\gamma_{-k},\beta_{-k},\alpha_{-k},\mu^{(r)},y_{1:T})\propto \prod_{t}{\mathcal{N}}(y_{t};{\xi_{t}^{(r)}}'\mu^{(r)},\sigma_{t}^{2})
\end{equation}
where $\gamma_{-k}=(\gamma_{1},\dots,\gamma_{k-1},\gamma_{k+1},\dots,\gamma_{M})$, $\beta_{-k}=(\beta_{1},\dots,\beta_{k-1},\beta_{k+1},\dots,\beta_{M})$ and $\alpha_{-k}=(\alpha_{1},\dots,\alpha_{k-1},\alpha_{k+1},\dots,\alpha_{M})$. 
We now follow the ARMA approximation of regime specific GARCH process i.e. 
\be
\begin{aligned}
\sigma_{t}^{2} &=\xi_{t}'\gamma + (\xi_{t}'\alpha)\epsilon_{t-1}^{2} + (\xi_{t}'\beta)\sigma_{t-1}^{2}\\
\epsilon_{t}^{2}&= \xi_{t}'\gamma + (\xi_{t}'\alpha+\xi_{t}'\beta)\epsilon_{t-1}^{2} - (\xi_{t}'\beta)(\epsilon_{t-1}^{2}-\sigma_{t-1}^{2}) + (\epsilon_{t}^{2}-\sigma_{t}^{2}).
\end{aligned}
\ee
Let
$$
w_{t} = \epsilon_{t}^{2}-\sigma_{t}^{2} = \left(\dfrac{\epsilon_{t}^{2}}{\sigma_{t}^{2}}-1\right)\sigma_{t}^{2} = (\chi^{2}(1)-1)\sigma_{t}^{2}
$$
with 
$$E_{t-1}[w_{t}]=0; \quad {\text{and}}\quad Var_{t-1}[w_{t}]=2\sigma_{t}^{4}.$$
Subject to the above and following \cite{Nak98} suggestion, we assume that $w_{t}\approx w_{t}^{*} \sim \mathcal{N}(0,2\sigma_{t}^{4})$. Then we have an \lq\lq auxiliary\rq\rq ARMA model for the squared error $\epsilon_{t}^{2}$.
\begin{equation}
\begin{aligned}
\epsilon_{t}^{2}&= \xi_{t}'\gamma + (\xi_{t}'\alpha+\xi_{t}'\beta)\epsilon_{t-1}^{2} - (\xi_{t}'\beta)w_{t-1}^{*} + w_{t}^{*}, \quad w_{t}^{*} \sim \mathcal{N}(0,2\sigma_{t}^{4})\\
{\text{i.e.}}\quad  w_{t}^{*} &= \epsilon_{t}^{2} - \xi_{t}'\gamma -  (\xi_{t}'\alpha)\epsilon_{t-1}^{2} - (\xi_{t}'\beta)(\epsilon_{t-1}^{2}-w_{t-1}^{*})
\end{aligned}
\end{equation}
Following \cite{Ardia08} we further express $w_{t}^{*}$ as a linear function of $(3\times1)$ vector of $(\gamma_{k},\alpha_{k},\beta_{k})'$. To do this, we approximate the function $w_{t}^{*}$ by first order Taylor's expansion about $(\gamma_{k}^{(r-1)},\alpha_{k}^{(r-1)},\beta_{k}^{(r-1)})'$.
\be
w_{t}^{*} \approx w_{t}^{**} =w_{t}^{*}(\theta_{-\pi}^{(r-1)}) - ((\gamma_{k},\alpha_{k},\beta_{k}) - (\gamma_{k}^{(r-1)},\alpha_{k}^{(r-1)},\beta_{k}^{(r-1)}))\nabla_{t} 
\ee
where
\be
\begin{aligned}
\dfrac{\partial w_{t}^{*}}{\partial\gamma_{k}} &= - \xi_{tk} + (\xi_{t}'\beta)\dfrac{\partial w_{t-1}^{*}}{\partial\gamma_{k}}\\
\dfrac{\partial w_{t}^{*}}{\partial\alpha_{k}} &= - \xi_{tk}\epsilon_{t-1}^{2} + (\xi_{t}'\beta)\dfrac{\partial w_{t-1}^{*}}{\partial\alpha_{k}}\\
\dfrac{\partial w_{t}^{*}}{\partial\beta_{k}}  &= - \xi_{tk}(\epsilon^{2}_{t-1}-w_{t-1}^{*}) + (\xi_{t}'\beta)\dfrac{\partial w_{t-1}^{*}}{\partial\beta_{k}}
\end{aligned}
\ee
\be
\nabla_{t} = -\left(\dfrac{\partial w_{t}^{*}}{\partial\gamma_{k}}, \dfrac{\partial w_{t}^{*}}{\partial\alpha_{k}},\dfrac{\partial w_{t}^{*}}{\partial\beta_{k}}\right)'\left|\right._{(\gamma_{k}=\gamma_{k}^{(r-1)},\alpha_{k}=\alpha_{k}^{(r-1)},\beta_{k}=\beta_{k}^{(r-1)})}.
\ee
Upon defining $r_{t}^{*}=w_{t}^{*}(\theta_{-\pi}^{(r-1)})+(\gamma_{k}^{(r-1)},\alpha_{k}^{(r-1)},\beta_{k}^{(r-1)})\nabla_{t}$, it turns out that\\ $w_{t}^{**} = r_{t}^{*} - (\gamma,\alpha,\beta)\nabla_{t}$. Furthermore, by defining the $T\times 1$ vectors\\ ${\bf{w}}=(w_{1}^{**},\dots,w_{T}^{**})'$, ${\bf{r^{*}}}=(r_{1}^{*},\dots,r_{T}^{*})'$ and ${\bf{\nabla}}=(\nabla_{1},\dots,\nabla_{T})'$ as well as a $T\times T$ matrix 
\be
{\bf{V}} = 2
\begin{pmatrix}
\sigma_{1}^{**4} &\cdots & 0\\
\vdots         &\ddots & \vdots\\
0              &\cdots & \sigma_{T}^{**4}\\
\end{pmatrix}
\ee
with $\sigma_{t}^{**2}=({\xi_{t}^{(r)}}'\gamma^{(r-1)})+({\xi_{t}^{(r)}}'\alpha^{(r-1)})(y_{t-1}-{\xi_{t-1}^{(r)}}'\mu^{(r)})^{2}+({\xi_{t}^{(r)}}'\beta^{(r-1)})\sigma_{t-1}^{**2}$,\\
we can approximate the full conditional probability of the regime specific volatility parameters as
\begin{equation}
\begin{aligned}
f(\gamma_{k},\beta_{k},\alpha_{k}|\xi_{1:T}^{(r)},\gamma_{-k},\beta_{-k},\alpha_{-k},\mu^{(r)},y_{1:T})&\propto \dfrac{1}{|{\bf{V}}|^{\frac{1}{2}}} \exp{\left(-\dfrac{{\bf{w'V^{-1}w}}}{2}\right)}\\
              &= {\mathcal{N}}_{3}(\mu,\Sigma)|_{\gamma_{k}>0,\alpha_{k}>0,\beta_{k}>0}
\end{aligned}\label{eqGARCHprop}
\end{equation}
where 
\be
\begin{aligned}
\Sigma &= (\nabla'{\bf{V}}^{-1}\nabla)^{-1}\\
\mu    &= \Sigma\nabla'{\bf{V}}^{-1}{\bf{r}}^{*}.
\end{aligned}
\ee
To sample for the truncated multivariate Normal distribution given in equation (\ref{eqGARCHprop}), we implement the Gibbs sampling technique by \cite{Ste12} for sampling from a truncated multivariate Normal distribution.
\clearpage
\bibliographystyle{plainnat}
\bibliography{MSGARCH}
\end{document}